\documentclass[12pt]{amsart}

\usepackage{amsmath}
\usepackage{amsthm}
\usepackage{latexsym}
\usepackage{amssymb}
\usepackage{mathrsfs}
\usepackage{tikz}

\newtheorem{theorem}{Theorem}
\newtheorem{proposition}{Proposition}

\newtheorem{definition}{Definition}
\newtheorem{example}{Example}
\newtheorem{remark}{Remark}

\makeatletter
\@namedef{subjclassname@2020}{\textup{2020} Mathematics Subject Classification}
\makeatother
\def\jpopn#1#2{%
  \mathopen{%
    \setbox0=\hbox{$#1\langle$}%
    \setbox2=\hbox{%
            {\hbox{$#1\langle$}}%
            \kern -.6\wd0\box0%
    }%
            \box2%
  }%
}

\def\jpcls#1#2{%
  \mathclose{%
    \setbox0=\hbox{$#1\rangle$}%
    \setbox2=\hbox{%
            {\hbox{$#1\rangle$}}%
            \kern -.6\wd0\box0%
    }%
            \box2%
  }%
}

\def\tp{\ {}^{t}\kern -3pt}

\allowdisplaybreaks
\AtBeginDocument{
  \renewcommand{\setminus}{\mathbin{\backslash}}%
}
\def\np#1{\ensuremath{\mathbf{NP}(\kern -1pt\it{#1})}}
\def\dnp#1{\ensuremath{\partial\mathbf{NP}(\kern -1pt\it{#1})}}

\renewcommand{\kappa}{\ensuremath{\varkappa}}
\renewcommand{\phi}{\ensuremath{\varphi}}
\renewcommand{\epsilon}{\ensuremath{\varepsilon}}

\begin{document}
%
\title[Microlocal regularity of the analytic vectors]{On the microlocal regularity of the analytic vectors for ``sums of squares" of vector fields}
\author{Gregorio Chinni}
\address{5, Via Val d'Aposa, 40123 Bologna, Italy}
\email[G. Chinni]{gregorio.chinni@gmail.com}
\thanks{The first author was partially supported by the Austrian Science Fund (FWF), Lise-Meitner position, project no. M2324-N35.} 
\author{Makhlouf Derridj}
\address{5, Rue de la Juvini\'ere, 78350 Les Loges en Josas, France}
\email[M. Derridj]{makhlouf.derridj@outlook.fr}
\begin{abstract} 
	We prove via FBI-transform a result concerning the microlocal Gevrey regularity 
	of analytic vectors for operators sums of squares of vector fields with real-valued real analytic coefficients of H\"ormander type,
	thus providing a microlocal version, in the analytic category, of a result due to M. Derridj in \cite{D_2019}
	concerning the problem of the local regularity for the Gevrey vectors for sums of squares of vector fields
	with real-valued real analytic/Gevrey coefficients.
	%

\vspace{0.5em}
	Nous d\'emontrons , en utilisant la transformation de Fourier-Bros-Iagolnitzer, 
	un r\'esultat de r\'egularit\'e Gevrey microlocale , optimale, des vecteurs analytiques d'op\'erateurs de H\"ormander de type
	"Sommes de carr\'es de champs de vecteurs" \`a coefficients analytiques sur un ouvert. 
	Ce r\'esultat est, dans le cadre analytique, la version microlocale du r\'esultat de M.Derridj \cite{D_2019},
	obtenu pour les vecteurs de Gevrey de tels op\'erateurs \`a coefficients Gevrey. 
\end{abstract}    
%
%
\maketitle
\section{Introduction}
\noindent
	We deal with the microlocal regularity of the analytic vectors for sum of squares of vector fields.
	Let $X_{1}(x,D), \, \dots, \, X_{m}(x,D)$ be vector fields with real-valued real analytic coefficients on $U$,
	open neighborhood of the origin in $\mathbb{R}^{n}$. Let $P(x,D)$ denote the corresponding sum of squares operator
%
\begin{align}\label{Op_P}
P(x,D) = \sum_{j=1}^{m}X_{j}^{2}(x,D).
\end{align}  
%
	We assume that the operator $P$ satisfies the H\"ormander's condition: the Lie algebra generated by the vector fields and
	their commutators has the dimension $n$, equal to the dimension of the ambient space.\\
	The operator $ P $ satisfies the \textit{a priori} estimate
%
%
\begin{equation} \label{eq:Cinftysubell}
\| u \|_{1/r}^{2} + \sum_{j=1}^{m} \|X_{j} u \|_{0}^{2} 
\leq C \left( |\langle P u , u \rangle | + \| u\|^{2}_{0} \right),
\end{equation}
%
	which we call, for the sake of brevity, the ``subelliptic estimate.''
	Here $ u \in C_{0}^{\infty}(U) $, $ \| \cdot \|_{0}  $ denotes the norm in $ L^{2}(U) $ and $ \| \cdot  \|_{s} $ the Sobolev norm of order $ s $ in $ U $. 
	Here  $ r $ is the least integer such that the vector fields, the commutators, the triple commutators etcetera up to the
	commutators of length $ r $ span at any point of the closure of $U$ all the ambient space $\mathbb{R}^{n}$.  
	The  sub-elliptic estimate was proved first by H\"ormander in \cite{H67} for a Sobolev norm of order $r^{-1}+\varepsilon$
	and up to order $r^{-1}$ subsequently by Rothschild and Stein \cite{RS} as well as in a pseudodifferential context
	by Bolley, Camus and Nourrigat in \cite{BCN-82}.\\

	\noindent
	Let $X_{j}(x,\xi)$ be the symbol of the vector field $X_{j}$.
	Write $\lbrace X_{i},X_{k}\rbrace$ the Poisson bracket of the symbols of the vector fields $X_{i}$, $X_{k}$:
%
\begin{align*}
\lbrace X_{i}, X_{k}\rbrace \left(x,\xi\right) =
\sum_{\ell=1}^{n} \left( \frac{\partial X_{i}}{\partial \xi_{\ell}} \frac{\partial X_{k}}{\partial x_{\ell}} - 
\frac{\partial X_{k}}{\partial \xi_{\ell}} \frac{\partial X_{i}}{\partial x_{\ell}}\right) (x,\xi).
\end{align*}
\begin{definition}\label{length}
	Let $\!(x_{0}, \xi_{0})\!$ be a point in the characteristic set of $ P $:
	\begin{align}
	\text{Char}(P)
	= \lbrace (x, \xi) \in T^{*}U\setminus \lbrace 0 \rbrace \, : \, X_{j}(x,\xi) = 0,\, j=1,\dots\, m \rbrace.
	\end{align}
	Consider all the iterated Poisson brackets $\lbrace X_{i}, X_{k} \rbrace $, $\lbrace X_{p}, \lbrace X_{i}, X_{k} \rbrace \rbrace $  etcetera.
	We define $\nu(x_{0},\xi_{0})$ as the length of the shortest iterated Poisson bracket of the symbols of the vector fields
	which is non zero in $(x_{0},\xi_{0})$.
\end{definition}
%
	We recall
%
\begin{definition}
	Let $P(x,D)$ be as in (\ref{Op_P}).
	We denote by $G^{s}(U; P)$ which is the space of the Gevrey vectors of order $s$ with respect to $P$,
	the set of all distributions $u \in \mathscr{D}'(U)$ such that for any compact subset $K$ of $U$
	there exists a positive constant $C_{K}$ such that
	\begin{align}\label{G_s-vectors_S_S}
	\|P^{N}u\|_{L^{2}(K)} \leq C_{K}^{2N+1} ((2N)!)^{s}, \qquad \forall\, N\in \mathbb{Z}_{+}.
	\end{align}
	When $s=1$ we set $G^{1}(U; P)= \mathscr{A}(U; P)$ the set of the analytic vectors with respect to $P$.
\end{definition}
%
	We recall that concerning systems of vector fields with real analytic coefficients satisfying H\"ormander's condition the problem of the local 
	regularity of the analytic vectors for such systems was first studied in  \cite{DH-1980} followed by a more refined version in \cite{HM_1980}.
	\\
	In a couple of recent works  M. Derridj, \cite{D_2018} and \cite{D_2019}, studied the problem of the local regularity for the Gevrey vectors
	for operators of H\"ormander type of first kind, i.e. sum of squares, and of the second kind or degenerate elliptic parabolic. 
	We prove  the minimal microlocal version of the result in \cite{D_2018} in the case of analytic vectors: 
%
%
\begin{theorem}\label{Miro_AnV}
	Let $P$ be as in (\ref{Op_P}). Let $u$ be an analytic vector for $P$, $u \in  \mathscr{A}(U; P)$. Let $(x_{0},\xi_{0})$ be a point in the characteristic
	set of $P$ and $\nu(x_{0},\xi_{0})$ its length. Then $(x_{0},\xi_{0}) \notin WF_{\nu(x_{0},\xi_{0})}(u)$.   
\end{theorem}
\noindent
	Where $ WF_{s}(u)$, $ s\geq 1$, denotes the wave front set of the distribution $u$; it will be defined in the next section via FBI-transform,
	Definition \ref{def:swf}. 
%
\begin{remark}
	A few remarks are in order:
	\begin{itemize}
		\item[i)] the method used to gain the above result can be extended to a class of H\"ormader type operators
		not strictly sums of squares; we consider operators of the form 
		$P(x,D) + \sum_{i=1}^{m}b_{j}(x)X_{j}(x,D) + c(x)$ where $P$ is as in (\ref{Op_P}), $b_{j}(x) $ are real-valued real analytic functions 
		and $c(x)$ is a real analytic complex function;
		\item[ii)] the strategy to obtain the above result can be carried over
		to the case of $s$-Gevrey vectors with $s \in \mathbb{Z}_{+}$;
		\item[iii)] the result is optimal, see example given in \cite{bccj_2016}.
	\end{itemize}
\end{remark}
%
	A few words about the method of proof: it consists
	in using the FBI transform and the subelliptic inequality on the FBI
	side obtained in \cite{ABC}. To do that we use a deformation technique of the
	Lagrangean associated to the FBI proposed by Grigis and Sj\"ostrand in \cite{GS}.
%

\vspace{1.8em}

\noindent
\textbf{Acknowledgement.} The authors would like to thank Antonio Bove for his comments and suggestions
in order to improve the manuscript. 	
%
\section{ Background on FBI and Micro-local sub-elliptic estimate for Sums of Squares}
%
	We are going to use a pseudodifferential and FIO (Fourier Integral
	Operators) calculus introduced by Grigis and Sj\"ostrand in the paper
	\cite{GS}. We recall below the main definitions
	and properties to make this paper self-consistent and readable. For
	further details we refer to the paper \cite{GS} and notes \cite{Sj-Ast}.
%
\bigskip
\paragraph{{\bf FBI Transform.}}
	We define the \textit{ FBI transform} of a temperate distribution $ u \in \mathscr{E}'(\Omega)$, $\Omega$ open subset
	of $\mathbb{R}^{n}$, as
%
\begin{align*}
Tu(z, \lambda) = \int_{\mathbb{R}^{n}} e^{i \lambda \psi(z, y)} u(y) dy,
\end{align*}
%
	where $ z \in \mathbb{C}^{n} $, $ \lambda \geq 1 $ is a large parameter, $ \psi(z,w) $ in $\mathbb{C}^{2n}$ is an holomorphic function
	such that $ \det \partial_{z}\partial_{w}\psi \neq 0 $, $ \Im \partial_{w}^{2}\psi > 0 $.
	To the phase $ \psi $ there corresponds a weight function $ \phi(z) $, defined as
%
\begin{align*}
\phi(z) = \sup_{y \in \mathbb{R}^{n}} - \Im \psi(z, y), \qquad z \in \mathbb{C}^{n}.
\end{align*}
\begin{example} \label{fbiclassic}
	A typical phase function may be $ \psi(z, y) = \frac{i}{2} (z-y)^{2} $. The corresponding weight fun\-ction is given by
	$ \phi(z)\doteq \phi_{0}(z) = \frac{1}{2} (\Im z)^{2} $.
\end{example}
%
	We recall that $T$ is associated to the following complex canonical transformation: 
%
\begin{align}\label{Can_Tr}
\begin{matrix}
\mathscr{H}_{T} \colon \mathbb{C}_{(w,\theta)}^{2n} \longrightarrow \mathbb{C}_{(z,\zeta)}^{2n},
\\
\quad \left(w, - \partial_{w} \psi(z,w)\right) \mapsto \left(z, \partial_{z} \psi(z,w)\right), 
\end{matrix}
\end{align}
%
	with $\psi$ as a generating function. \\
	In particular $\mathscr{H}_{T}(\mathbb{R}^{2n})\doteq \Lambda_{\phi}= \lbrace \left( z,-2i \partial_{z}\phi(z)\right); z\in \mathbb{C}^{n}\rbrace$. 
	In the case of classical phase function, see Example \ref{fbiclassic}, we have
%
\begin{align*}
\mathscr{H}_{0} (x, \xi)= ( x-i\xi, \xi ), \qquad (x,\xi) \in \mathbb{R}^{2n}.
\end{align*}
%
	We set $\mathscr{H}_{0} (\mathbb{R}^{2n}) = \Lambda_{\phi_{0}}$.\\
	We recall the definition of $ s $--Gevrey wave front set of a distribution via classical FBI transform,
	i.e. using the phase function and the corresponding weight function of the Example \ref{fbiclassic}.
%
\begin{definition}
	\label{def:swf}
	Let $u$ be a compactly supported distribution on $\mathbb{R}^{n}$.
	Let $ (x_{0}, \xi_{0}) \in  T^{*}\mathbb{R}^{n}\setminus 0 $. We say
	that $ (x_{0}, \xi_{0}) \notin WF_{s}(u) $, $s \geq 1$, if there exist a
	neighborhood $ \Omega $ of $ x_{0} - i \xi_{0} \in \mathbb{C}^{n} $ and
	positive constants $ C $, $\varepsilon $ such that
	$$ 
	| e^{-\lambda \phi_{0}(z)} Tu(z, \lambda) | \leq C e^{-\varepsilon \lambda^{1/s}},
	$$
	for every $ z \in \Omega $ and $\lambda > 1$.
\end{definition}
\bigskip
%
\paragraph{{\bf Pseudodifferential Operators.}}
%
	Let us consider $ (z_{0}, \zeta_{0})$  $\in \mathbb{C}^{2n} $ and a real valued real analytic
	function $ \phi(z) $ defined near $ z_{0} $, such that $ \phi $ is strictly plurisubharmonic and
%
\begin{align*}
\frac{2}{i}\ \partial_{z}\phi(z_{0}) =\zeta_{0}.
\end{align*}
%
	Denote by $ \vartheta(z, w) $ the holomorphic function defined near $ (z_{0}, \bar{z}_{0}) $ by 
%
\begin{align}\label{psi}
\vartheta(z, \bar{z}) = \phi(z).
\end{align}
%
	Because of the strict plurisubharmonicity of $ \phi $, we have 
%
\begin{align}\label{psixy}
\det \partial_{z} \partial_{w} \vartheta \neq 0 
\end{align}
%
	and
%
\begin{align}
\label{repsi}
\Re \vartheta(z, \bar{w}) - \frac{1}{2}\left [ \phi(z) + \phi(w) \right] \sim - |z - w |^{2}.
\end{align}
%
	Let $ \lambda \geq 1 $ be a large positive parameter. We write 
%
\begin{align*}
\tilde{D} = \frac{1}{\lambda}  D, \qquad  D = \frac{1}{i} \partial.
\end{align*}
%
	Denote by $ q(z, \zeta, \lambda) $ an analytic classical symbol and by $ Q(z, \tilde{D}, \lambda) $ the formal classical
	pseudodifferential operator associated to $ q $.
	Using ``Kuranishi's trick''
	\footnote{For more details on the ``Kuranishi's trick'' see \cite{H_FIO}
		Proposition 2.1.3 and \cite{Sj-Ast} Remarque 4.3.}
	one may represent  $ Q(z, \tilde{D},\lambda) $ as
	\begin{equation}\label{qk}
	Q u(z, \lambda) = \left( \frac{\lambda}{2 i\pi }\right )^{n} 
	\int e^{2\lambda (\vartheta(z, \theta) - \vartheta(w, \theta))} \tilde{q}(z, \theta, \lambda) u(w) dw d\theta.
	\end{equation}
	Here $ \tilde{q} $ denotes the symbol of $ Q $ in the actual representation. 
	
	To realize the above operator we need a prescription for the
	integration path\footnote{For a detailed discussion about the integration paths see \cite{Sj-Ast}.}.
	This is accomplished by transforming the classical integration path
	via the Kuranishi change of variables and eventually applying Stokes
	theorem: 
	\begin{equation}
	\label{realization} 
	Q^{\Omega}u(z, \lambda) = \left( \frac{\lambda}{\pi}\right)^{n}
	\int_{\Omega} e^{2 \lambda \vartheta(z, \bar{w})} \tilde{q}(z, \bar{w}, \lambda) u(w) e^{-2 \lambda\phi(w)} L(dw),
	\end{equation}
	where $ L(dw) = (2i)^{-n} dw \wedge d\bar{w} $ is the Lebesgue measure in $\mathbb{R}^{2n}$,
	the integration path is $ \theta = \bar{w} $ and $ \Omega \times \overline{\Omega} $
	is a small neighborhood of $ (z_{0}, \bar{z}_{0}) $. We remark that $Q^{\Omega}u (z)$ is an holomorphic
	function of $z$. 
%
\begin{definition}	\label{def:Hphi}
	Let $ \Omega $ be an open subset of $ \mathbb{C}^{n} $.
	We denote by $ H_{\phi}(\Omega) $ the space of all functions $ u(z, \lambda) $ holomorphic with respect to $z$, such that for every $ \epsilon > 0 $
	and for every compact $ K \subset\!\subset \Omega $ there exists a constant $ C > 0 $ such that
	\begin{align*}
	| u(z, \lambda) | \leq C e^{\lambda (\phi(z) + \varepsilon)},
	\end{align*}
	for $ z \in K $ and $ \lambda \geq 1 $.
\end{definition}
%
	A few remarks are in order.
	
	\vspace*{0.2em}
	\begin{enumerate}
		\item[i)] \label{rem:1}
		If $ \tilde{q} $ is a classical symbol of order zero, $ Q^{\Omega} (z, \tilde{D},\lambda)$
		is uniformly bounded as $ \lambda \rightarrow +\infty $, from $ H_{\phi}(\Omega) $ into itself.
		\vspace*{0.2em}
		\item[ii)] \label{rem:2}
		If the principal symbol is real, $ Q^{\Omega} (z, \tilde{D},\lambda)$ is formally self adjoint operator in $L^{2}(\Omega,$  $ e^{ -2 \lambda \phi(z)} L(dz))$.
		\vspace{0.2em}
		\item[iii)] \label{rem:3}
		The definition \eqref{qk} of the realization of a pseudodifferential
		operator on an open subset $ \Omega $ of $ \mathbb{C}^{n} $ is not the
		classical one. Via the Kuranishi trick it can be reduced to the
		classical definition. On the other hand using the function $ \vartheta $ allows us to use a weight function not explicitly related to an FBI
		phase. This is useful since in the proof we deform the I-Lagrangian, R-Symplectic variety $\Lambda_{\phi_{0}} $, corresponding e.g. to the classical FBI phase,
		and obtain a \textit{deformed} weight function which is useful in the a priori estimate.
	\end{enumerate}
	%
	We also recall that the identity operator can be realized as
%
\begin{align}
\label{idomega}
I^{\Omega}u(z, \lambda) = \left( \frac{\lambda}{\pi}\right)^{n}
\int_{\Omega} e^{2 \lambda \vartheta(z, \bar{w})} i(z, \bar{w}, \lambda) u (w, \lambda) e^{-2\lambda \phi(w)} L(dw),
\end{align}
%
	for a suitable analytic classical symbol $ i(z, \zeta, \lambda) $. Moreover we have the following estimate (see \cite{GS} and \cite{Sj-Ast})
%
\begin{align}
\label{errorest}
\| I^{\Omega} u - u \|_{\phi - d^{2}/C} \leq C' \| u \|_{\phi + d^{2}/C},
\end{align}
%
	for suitable positive constants $ C $ and $ C' $. Here we denoted by
%
\begin{align}\label{d}
d(z) = \text{dist}(z, \complement \Omega),
\end{align}
%
	the distance of $ z $ to the boundary of $ \Omega $, and by
%
\begin{align}
\label{eq:phinorm}
\| u \|_{\phi}^{2} = \int_{\Omega} e^{-2\lambda \phi(z)} |u(z) |^{2} L(dz).
\end{align}
%
	We also recall the following important result on the composition of two pseudodifferential operators.
%
\begin{proposition}[\cite{GS}]\label{composition}
	Let $ Q_{1} $ and $ Q_{2} $ be of order zero. Then they can be composed and
	$$ 
	Q_{1}^{\Omega} \circ Q_{2}^{\Omega} = (Q_{1} \circ Q_{2})^{\Omega} + R^{\Omega},
	$$
	where $ R^{\Omega} $ is an error term, i.e. an operator whose norm is
	$ \mathscr{O}(1) $ as an operator from $ H_{\phi + (1/C) d^{2}} $ to $	H_{\phi - (1/C) d^{2}} $ 
\end{proposition}
%
\bigskip
\paragraph{{\bf The \textit{a priori} Estimate}}
%
	Let $X_{j}(z,\zeta)$, $j=1, \, \dots, \, m $, be classical analytic symbols of order one defined in $\Omega$
	open neighborhood of $(z_{0},\zeta_{0}) \in \Lambda_{\phi}$ in $\mathbb{C}^{2n}$. We assume also that the
	$X_{j|_{\Lambda_{\phi}} }$ are real valued. Let 
%
\begin{align}
P(z,\tilde{D}) = \sum_{j=1}^{m} X_{j}^{2}(z,\tilde{D}).
\end{align}
%
	According to \cite{GS} the $ \Omega $-realization of $ P $  can be written as
%
\begin{align}\label{eq:Preal}
P^{\Omega} = \sum_{j=1}^{m} (X_{j}^{\Omega})^{2} + \mathscr{O}(\lambda^{2}),
\end{align}
%
	where $\mathscr{O}(\lambda^{2}) $ is continuous from $
	H_{\tilde{\phi}} $ to $ H_{\phi - (1/C)d^{2}} $ with norm boun\-ded by $
	C' \lambda^{2} $, $ \tilde{\phi} $ given by
%
\begin{align*}
\tilde{\phi}(z) = \phi(z) + \frac{1}{C} d^{2}(z),
\end{align*}
%
	and $ d $ has been defined in (\ref{d}).\\
	Following \cite{ABC} we state the FBI version of the estimate \eqref{eq:Cinftysubell}.
%
\begin{theorem}	\label{Micro_Sub_Est} 
	Let $(x_{0},\xi_{0}) $ be in $\text{Char}(P)$ and $\nu \doteq \nu(x_{0},\xi_{0})$, Definition \ref{length}.
	Let $\mathscr{H}_{T}(x_{0},\xi_{0})=(z_{0},\zeta_{0}) \in \Lambda_{\phi}$ and $P^{\Omega}$ be as in (\ref{eq:Preal}).
	Let $ \Omega_{1}   $ open neighborhood of $(z_{0}, \zeta_{0})$ such that $\Omega_{1}\subset\!\subset \Omega $.
	Then
	\begin{equation}
	\label{apriorifin}
	\lambda^{\frac{2}{\nu}} \| u \|_{\phi}^{2} +
	\sum_{j=1}^{m} \| X_{j}^{\Omega} u\|_{\phi}^{2} 
	\leq C \left( \langle P^{\Omega}u , u\rangle_{\phi} + \lambda^{\alpha} \|u\|_{\phi, \Omega\setminus\Omega_{1}}^{2} \right),
	\end{equation}
	where $ \alpha $ is a positive integer and $ u \in L^{2}(\Omega, e^{-2 \phi(z)}L(dz)) $. 
\end{theorem}
%
%
\section{Proof of the Theorem \ref{Miro_AnV}}
%
	In order to prove the result we want take advantage of Theorem \ref{Micro_Sub_Est}.
	We consider the sum of squares operator
%
\begin{align}\label{Op_Q}
Q(x, D_{t}, D)= \sum_{j=0}^{m} X_{j}^{2} = D_{t}^{2}+P(x,D),
\end{align}
%
	in $\tilde{\mathcal{O}}=\left] -\delta_{0},\delta_{0} \right[\times \mathcal{O}$, $\delta_{0} > 0$.\\
	We study the microlocal properties of the solutions of the problem $Qv= f$, $f \in C^{\omega}(\tilde{\mathcal{O}})$.  
	We denote by $\tilde{\Sigma}$ the characteristic set of $Q$ given by
%
\begin{align}
\tilde{\Sigma}
&= \lbrace (t,x,\tau,\xi) \in T^{*}\tilde{\mathcal{O}}\setminus\lbrace 0\rbrace \,:\, Q(t,x,\tau,\xi)=0\rbrace\\
\nonumber
&= \lbrace (t,x,\tau,\xi) \in T^{*}\tilde{\mathcal{O}}\setminus\lbrace 0\rbrace \,:\, \tau=0,\, X_{j}(x,\xi)=0,\, j=1,\, \dots, \, m\rbrace.
\end{align}
%
	We remark that $\nu_{(t_{0}, x_{0}, 0, \xi_{0})}$,$ (t_{0}, x_{0}, 0, \xi_{0}) \in \tilde{\Sigma}$, is equal to $ \nu_{(x_{0},\xi_{0})}$,  $(x_{0},\xi_{0}) \in \Sigma$,
	where $\Sigma $ denotes the characteristic set of $P(x,D)$.\\
	We construct a deformation of $\Lambda_{\phi_{0}}$ following the ideas in \cite{GS}, see also \cite{ABC}.
	Let $(0,x_{0}, 0,\xi_{0}) \in \tilde{\Sigma}$ and $\nu$ its length.\\ 
	We perform an FBI-transform of the form
%
\begin{align*}
Tu(z, \lambda) = \int_{\mathbb{R}^{n+1}} e^{ i\lambda \psi(z,t ,x )} u(t,x) dtdx,\qquad z=(z_{0},z_{1}) \in \mathbb{C}^{1+n},
\end{align*}
%
	where $ u(t,x) $ is a compactly supported distribution and $ \psi(z,t, x) $ is a phase function.
	Even though it does not really matter which phase function we use, the classical phase function will be employed:
%
\begin{align}\label{phaseFBI}
\psi_{0}(z, t, x) = \frac{i}{2}\left[ (z_{0}-t)^{2}+(z_{1}-x)^{2}\right]. 
\end{align}
%
	Let $ \Omega$ be an open neighborhood of the point $ \pi_{z}\circ\mathscr{H}_{T}(0,x_{0},0,\xi_{0}) $
	in $ \mathbb{C}^{1+n} $. Here $ \pi_{z} $ denotes the space projection
	$ \pi_{z} \colon  \mathbb{C}^{1+n}_{z} \times \mathbb{C}^{1+n}_{\zeta} \rightarrow \mathbb{C}^{1+n}_{z} $,
	$\zeta = (\zeta_{0},\zeta_{1})$, and $ \mathscr{H}_{T} $ is the complex canonical transformation associated to $ T $, (\ref{Can_Tr}).
	We recall that in the case of FBI with classical phase function, Example \ref{fbiclassic}, we have $ \mathscr{H}_{0} (t,x,\tau, \xi)= (t-i\tau, x-i\xi,\tau, \xi). $\\
	Denoting by $\tilde{Q}$ our operator after the FBI we have that $\tilde{Q}_{|\Lambda_{\phi_{0}}}= Q$,
	$\Lambda_{\phi_{0}} = \mathscr{H}_{0}(\mathbb{R}^{2(1+n)})$.
	We have that $\pi_{z} \circ \mathscr{H}_{0}(0,x_{0},0,\xi_{0})= (0,x_{0}-i\xi_{0}) = (0,w_{0}) \in  \mathbb{C}^{1+n}$.
	We perturb canonically $\phi_{0}$. For $ \lambda \geq 1 $ let us consider a real analytic function defined near the point
	$ \mathscr{H}_{0}(0,x_{0},0,\xi_{0}) \in \Lambda_{\phi_{0}} $, say $ h(z, \zeta, \lambda) $. 
	Solve, for small positive $ s $, the Hamilton-Jacobi problem
%
\begin{align}\label{eq:hj}
\left\{
\begin{matrix}
2 \dfrac{\partial \phi}{\partial s}(s, z, \lambda) =
h \left(z, \dfrac{2}{i} \dfrac{\partial \phi}{\partial z}(s, z, \lambda), \lambda \right) \\
\phi(0, z, \lambda) = \phi_{0}(z) \hfill
\end{matrix}
\right . .
\end{align}
%
	Set
%
\begin{align*}
\phi_{s}(z, \lambda) = \phi(s, z, \lambda),
\end{align*}
%
	we have the canonical map $\Lambda_{\phi_{0}} \rightarrow \Lambda_{\phi_{s}}$ where
%
\begin{align*}
\Lambda_{\phi_{s}} = \exp\left(i s H_{h}\right) \Lambda_{\phi_{0}}.
\end{align*}
%
	We choose the function $ h $ as
%
\begin{align*}
h (z, \zeta, \lambda) = h \left(z, \dfrac{2}{i} \dfrac{\partial \phi_{0}}{\partial z}( z), \lambda \right)
+ \lambda^{-1} h_{1}(z,\zeta) \left(\zeta - \dfrac{2}{i} \dfrac{\partial \phi_{0}}{\partial z}( z)\right),
\end{align*}
%
	where $h_{1}(z,\zeta)$ is an holomorphic function and 
%
\begin{align}\label{eq:h}
h \left(z, \dfrac{2}{i} \dfrac{\partial \phi_{0}}{\partial z}( z), \lambda \right) =
h (z, \zeta, \lambda)_{|\Lambda_{\phi_{0}}} = (z_{0}'')^{2}+\lambda^{-\frac{\nu-1}{\nu}} \left( (z_{0}')^{2}+ | z_{1} - w_{0} |^{2}\right), 
\end{align}
%
	$z_{0} = z_{0}'+iz_{0}'' \in \mathbb{C}$.\\
	Since $\mathbb{R}^{2(1+n)}$ and $\Lambda_{\phi_{0}}$ are isometric, keep in mind the definition of $ \Lambda_{\phi_{0}} $, 
	it is easier to construct the function $h$ in $ \mathbb{R}^{2(n+1)} $ near the characteristic point: 
%
\begin{align}\label{eq:hLphi0}
h (t, x, \tau, \xi, \lambda) = \tau^{2}+ \lambda^{-1+\frac{1}{\nu}} \left[ t^{2} + | x - x_{0} |^{2} + | \xi - \xi_{0} |^{2} \right].
\end{align}
%
	The function $ \phi_{s} $ can be expanded as a power series in the variable $ s $ using both equation \eqref{eq:hj}
	and the Fa\`a di Bruno formula to obtain
%
\begin{align}\label{eq:phi_s}
\phi_{s}(z, \lambda) = \phi_{0}(z) + \frac{s}{2}  \ h(\cdot, \cdot,
\lambda)_{\big |_{\Lambda_{\phi_{0}}}} + \mathscr{O}(\lambda^{-1}s^{2}) ,
\end{align}
%
	where $ h $ on $ \Lambda_{\phi_{0}} $ is given by \eqref{eq:h}.
	Our purpose is to use the estimate \eqref{apriorifin} where the weight
	function $ \phi $ has been replaced by the weight $ \phi_{s} $. This
	is possible using the phase $ \vartheta_{s} $ in \eqref{qk} and realizing
	the operator as in \eqref{realization}. Here $ \vartheta_{s} $ is defined
	as the holomorphic extension of $ \vartheta_{s}(z, \bar{z}) = \phi_{s}(z)$. 
	
	We need to restrict the symbol of $ Q $ to $ \Lambda_{\phi_{s}} $; we denote by $ Q^{s} $ the symbols of $ Q $ restricted to $\Lambda_{\phi_{s}} $.
	Noting that
%
\begin{multline*}
X_{j}^{2}(x, \frac{2}{i} \partial_{x}\phi_{s}(x, \lambda), \lambda) =
X_{j}^{2}(x, \frac{2}{i} \partial_{x}\phi_{0}(x), \lambda) 
\\
\qquad+ 2 s
X_{j}(x, \frac{2}{i} \partial_{x}\phi_{0}(x), \lambda)
\langle \partial_{\xi} X_{j}(x,
\frac{2}{i} \partial_{x}\phi_{0}(x),
\lambda), \frac{2}{i} \partial_{x}\partial_{s} {\phi_{s}(x, \lambda)}_{\big|_{s=0}} \rangle 
\\
+ \mathscr{O}(s^{2} \lambda^{\frac{2}{\nu}}). 
\end{multline*}
%
	We deduce that
%
\begin{align}\label{lt}
Q^{s}
&=Q+ s\sum_{j=0}^{m} X_{j}\lbrace h, X_{j}\rbrace +s^{2}\sum_{j=0}^{m}\lbrace h, X_{j}\rbrace^{2}
+\mathcal{O}(s^{2}\lambda^{\frac{2}{\nu}})
\\
&
\nonumber
\qquad
=  Q(x, \xi) + s R(x, \xi, \lambda) +
\mathscr{O}(s^{2} \lambda^{\frac{2}{\nu}}).
\end{align}
%
	The analytic extension of $ Q^{s} $ is the symbol appearing in the $\Omega $-realization of $ Q^{s} $, $ {Q^{s}}^{\Omega} $. We point out
	that the principal symbol of $ Q^{s} $ satisfies the assumptions of Theorem \ref{Micro_Sub_Est} and, using the a priori inequality
	(\ref{apriorifin}), we can deduce an estimate of the form (\ref{apriorifin}) for $ Q^{s} $ in the $ H_{\phi_{s}} $ spaces.
	We have
%
\begin{align*}
&\lambda^{\frac{2}{\nu}}\| u \|_{\phi_{s}}^{2} + \sum_{j=0}^{m}
\|X_{j}^{\Omega} u\|^{2}_{\phi_{s}} 
\\
&
\hspace*{3em}
\leq C \left( |\langle ({Q^{s}}^{\Omega} - s R^{\Omega} - \mathscr{O}(s^{2} \lambda^{\frac{2}{\nu}}))u
, u\rangle_{\phi_{s}}| + \lambda^{\alpha}
\|u\|_{\phi_{s}, \Omega\setminus\Omega_{1}}^{2} \right).
\end{align*}
%
	The third term in the right hand side of the scalar product above is
	easily absorbed on the left provided $ s $ is small enough. 
	Let us consider the second term in the scalar product above. 
	By Proposition \ref{composition}, we have
%
\begin{align*}
R^{\Omega} = \sum_{j=0}^{m} a_{j}^{\Omega}(x, \tilde{D}, \lambda)
X_{j}^{\Omega}(x, \tilde{D}, \lambda) + \mathscr{O}(\lambda),
\end{align*}
%
	where $ \mathscr{O}(\lambda) $ denotes an operator from $ H_{\phi_{s}+ \frac{1}{C} d^{2}} $ to $ H_{\phi_{s} - \frac{1}{C} d^{2}} $ whose
	norm is bounded by $ C \lambda $. Hence
%
\begin{align*}
s | \langle R^{\Omega}u, u\rangle_{\phi_{s}} | 
\leq C s 
\Big( \lambda^{\frac{2}{\nu}}\| u \|_{\phi_{s}}^{2} + \sum_{j=0}^{m} \|X_{j}^{\Omega} u\|^{2}_{\phi_{s}}
+ \lambda^{2} \|u \|^{2}_{\tilde{\phi}_{s}} \Big).
\end{align*}
%
	Hence we deduce that there exist a neighborhood $ \Omega_{0} $ of $ (0,w_{0})$, a positive
	number $ \delta $ and a positive integer $ \alpha $
	such that, for every $ \Omega_{1} \subset\!\subset \Omega_{2}
	\subset\!\subset\Omega \subset \Omega_{0} $, there exists a constant $ C > 0
	$ such that, for $ 0 < s < \delta $, we have 
%
\begin{align}
\label{scalprodt}
\lambda^{\frac{2}{\nu}}\| u \|_{\phi_{s}, \Omega_{1}} 
\leq C \left( \| {Q^{s}}^{\Omega} u\|_{\phi_{s}, \Omega_{2}} +
\lambda^{\alpha} \|u\|_{\phi_{s}, \Omega\setminus\Omega_{1}} \right). 
\end{align}
%
	We now prove that if $Qu$ is analytic at $(0,x_{0}, 0,\xi_{0})$ then the point $(0,x_{0}, 0,\xi_{0})$
	does not belong to $WF_{\nu}(u)$.
	Since $Qu$ is real analytic the first term in the right hand side of (\ref{scalprodt}) can be estimated
	by $C e^{-\lambda/C}$ for a positive constant $C$. We have to estimate the second term on the right hand side of the above inequality.
	We have
%
\begin{align*}
\phi_{s}(z,\lambda)=\phi_{0}(z) + \frac{s}{2}h(z,\frac{2}{i}\frac{\partial \phi}{\partial z}(0,z),\lambda) + \mathcal{O}(\lambda^{-1} s^{2}).
\end{align*}
%
	Hence
%
\begin{align*}
\phi_{s}(z,\lambda)-\phi_{0}(z) \sim \frac{s}{2} \left[ (z_{0}'')^{2}+\lambda^{-\frac{\nu-1}{\nu}} \left( (z_{0}')^{2}+ | z_{1} - w_{0} |^{2}\right) \right].
\end{align*}
%
	Since $ z= (z_{0},z_{1}) \in \Omega\setminus\Omega_{1}$, i.e. far from $(0,w_{0})$, there exists a positive constant
	$\beta$  such that 
%
\begin{align*}
h_{|\Lambda_{\phi_{0}} \cap  \Omega\setminus\Omega_{1}} \geq  2 \lambda^{-1+1/\nu}\beta > 0. 
\end{align*}
%
	We have
%
\begin{align*}
\phi_{s}(z,\lambda)_{|\Omega\setminus\Omega_{1}}\geq \phi_{0}(z) + s \lambda^{-1+\frac{1	}{\nu}}\beta +  \mathcal{O}(\lambda^{-1} s^2).
\vspace{0.3em}
\end{align*}
%
%
	The second term on the right hand side of (\ref{scalprodt}) can be estimated by
%
\begin{align*}
\|u\|_{\phi_{s},\Omega\setminus\Omega_{1}}^{2}\leq C_{1}(s) e^{-\lambda^{ \frac{1}{\nu}} s\frac{\beta}{2}},
\end{align*}
%
	where $C_{1}(s) > 0$. From (\ref{scalprodt}) and the above argument there is a positive constant $C_{2}$ such that
%
\begin{align*}
\|u\|_{\phi_{s},\Omega_{1}}^{2}\leq C_{2} 
e^{- \lambda^{\frac{1}{\nu}}s \frac{\beta}{2}}.
\end{align*}
%
	Let $\Omega_{3}$ a sufficient small neighborhood of the point $(0,w_{0})$,  $ \Omega_{3} \Subset \Omega_{1}$,
	such that for a fixed small positive $s$
%
\begin{align*}
\phi_{s}(z,\lambda)  - \phi_{0}(z) \leq  \frac{s\beta}{4} \lambda^{-1+\frac{1}{\nu}} + \lambda^{-1} C_{3}(s),
\end{align*}
$ z \in \Omega_{3} $, $ \lambda \geq 1 $.\\

\noindent
	Then there are two positive constants, $\tilde{C}$ and $\epsilon$, such that
%
\begin{align*}
\|u\|_{\phi_{0},\Omega_{3}}^{2}\leq \tilde{C} e^{-\epsilon \lambda^{\frac{1}{\nu}}}.
\end{align*}
\bigskip
%

	Now we consider the problem 
%
\begin{align}\label{eq:Q=0}
\left\{
\begin{matrix}
\left(D_{t}^{2} + P(x,D)\right) U(t,x)=0,\\
\\
U(0,x) =u(x),\qquad\qquad\quad
\end{matrix}
\right . 
\end{align}
%
	in $\tilde{\mathscr{O}}=\left] -\delta_{0},\delta_{0} \right[\times \mathscr{O}$, $\delta_{0} > 0$, where $u(x)$ is an analytic vector for $P(x,D)$:
%
\begin{align}\label{An_v_est}
\|P^{k}u\|_{0}\leq C^{2k+1} (2k)!.
\end{align}
%
	The function
%
\begin{align*}
U(t,x) = \sum_{k\geq 0} \frac{t^{2k}}{2k!} \, P^{k}u(x)
\end{align*}
%
	is a solution of the above problem. We choose $\delta_{0} < \sqrt{2}C$.\\
	In order to complete the proof of the Theorem \ref{Miro_AnV} we have to show that $(0,x_{0},0,\xi_{0})\notin WF_{s_{0}}\left(U\right)$
	if and only if $(x_{0},\xi_{0})\notin WF_{s_{0}}(u)$ for every $s_{0}\geq 1$.\\
	This result was showed in the case $s_{0}=1$ via Fourier transform in \cite{bcm_1990}, Proposition 3.3.
	We give, for any $ s_{0}  \in [1, +\infty )$, a proof via the classical FBI transform.\\ 
	\textbf{Step one}: if  $(x_{0},\xi_{0})\notin WF_{s_{0}}(u)$ then $(0,x_{0},0,\xi_{0}) \notin WF_{s_{0}}(U)$.\\
	By hypothesis we have that $(x_{0},\xi_{0}) \notin WF_{s_{0}}(u)$ if and only of
	there exist $\Omega$ open neighborhood of the point $x_{0}-i\xi_{0}$ in $\mathbb{C}^{n}$
	and positive constants $C_{1}$ and $\varepsilon_{1}$ such that
%
\begin{align} 
\label{FBI_u}
| e^{-\lambda \phi_{0}(z)} T\left(\chi u\right)(z, \lambda) | \leq C_{1} e^{- \varepsilon_{1}\lambda^{1/s_{0}}}, \qquad \forall z \in \Omega,
\end{align}
%
	where $\chi$ is a $C_{0}^{\infty}\left(\mathscr{O}\right)$ identically one in a neighborhood of $x_{0}$.\\
	We have to show that there is $\Xi$ open neighborhood of the point $(0,x_{0}-i\xi_{0})$ in $\mathbb{C}^{n+1}$
	and positive constants $C_{2}$ and $\varepsilon_{2}$ such that
	\begin{equation} 
	| e^{-\lambda \phi_{0}(w,z)} T\left(\chi U\right)(w,z, \lambda) | \leq C_{2} e^{- \varepsilon_{2}\lambda^{1/s_{0}}}, \qquad \forall (w,z) \in \Xi,
	\end{equation}
	where $\chi(t,x)= \chi_{_{0}}(t)\theta_{_{0}}(x)$, here $\chi_{_{0}}(t)$ is $C_{0}^{\infty}(]-\delta_{1}, \delta_{1}[)$, $0< \delta_{1}< \delta_{0}$,
	such that $\chi_{_{0}}(t)\equiv 1 $ in $]-\delta_{2}, \delta_{2}[$, $0< \delta_{2} < \delta_{1}/2$, and 
	$\theta_{_{0}}(x)$ is $C_{0}^{\infty}(\mathscr{B}_{r_{0}}(x_{0}))$, $\mathscr{B}_{r_{0}}(x_{0})= \lbrace x\in \mathbb{R}^{n}: | x-x_{0}| < r_{0}\rbrace$,
	$r_{0} \leq \text{dist}\left( x_{0}, \complement\pi_{z'}(\Omega)\right) $
	such that $\theta_{_{0}}(x)\equiv 1 $ in $\mathscr{B}_{r_{1}}(x_{0})$, $0< r_{1} < r_{0}$.\\
	We have
%
\begin{align*}
T\left(\chi U\right)(w,z,\lambda) =\iint e^{-\frac{\lambda}{2}(w-s)^{2}} e^{-\frac{\lambda}{2}(z-y)^{2}} \chi_{_{0}}(s)\theta_{_{0}}(y)U(s,y) \,dsdy
\\
= \sum_{N=0}^{\infty}\underbrace{ \frac{1}{(2N)!}\iint e^{-\frac{\lambda}{2}(w-s)^{2}} e^{-\frac{\lambda}{2}(z-y)^{2}} \chi_{_{0}}(s)\theta_{_{0}}(y)s^{2N}P^{N}u(y)
	\,dsdy}_{ \doteq \mathscr{P}_{N}(w,z)} \,.
\end{align*}
%
	Let $\tilde{r}_{0} > 0$ such that $\tilde{r}_{0} \ll r_{1}$. We take $z$ in the FBI transform such that
	$z'= \Re(z) \in \mathscr{B}_{\tilde{r}_{0}-\varepsilon}\left(x_{0}\right)$, where $ 0 < \varepsilon < \tilde{r}_{0}$.\\  
	Case $N=0$, since there are two positive constants $A $ and $\tilde{\varepsilon}_{0}$
	such that
%
\begin{align*}
\left|
\int e^{-\frac{\lambda}{2}(w-s)^{2}}\! \chi_{0}(s) ds 
\right|
\leq e^{\frac{\lambda}{2} \left( w''\right)^{2}} A e^{-\lambda\tilde{\varepsilon}_{0}},
\end{align*}
%
	taking advantage from (\ref{FBI_u}) there is a positive constant $C_{3}$ such that
%
\begin{align}\label{N=0}
&\left|e^{-\lambda \phi_{0}(w,z)}  \iint e^{-\frac{\lambda}{2}(w-s)^{2}} e^{-\frac{\lambda}{2}(z-y)^{2}} \chi_{_{0}}(s)\theta_{_{0}}(y)u(y) \,dsdy\right|
\\
\nonumber
&
\hspace{13em}
\leq  C_{3} e^{- \varepsilon_{1}\lambda^{1/s_{0}}}e^{- \tilde{\varepsilon}_{0}\lambda}.
\end{align}
%
	In order to make  the proof more readable, before looking at the general case, we analyze the cases $N=1$ and $N=2$.\\
	Case $ N=1 $; we have
%
\begin{align}\label{N=1}
\mathscr{P}_{1}\left(w,z\right)\doteq\frac{1}{2}\iint e^{-\frac{\lambda}{2}(w-s)^{2}} e^{-\frac{\lambda}{2}(z-y)^{2}} \chi_{_{0}}(s)\theta_{_{0}}(y)s^{2}Pu(y) \,dsdy.
\end{align}
%
	We introduce $\theta_{_{1}}(y)\in C_{0}^{\infty} \left(\mathscr{B}_{r_{1}}\left(x_{0}\right)\right)$
	such that $\text{supp}\left(\theta_{_{1}}\right)\subseteq \mathscr{B}_{r_{1}}\left(x_{0}\right)$, where $\theta_{_{0}}(y) \equiv 1$, and
	$\theta_{_{1}}(y) \equiv 1$ in $\mathscr{B}_{r_{2}}\left(x_{0}\right)$, where $\tilde{r}_{0} \leq r_{2}< r_{1}<r_{0}$. We have
%
\begin{align}\label{N=1_1}
\mathscr{P}_{1}\left(w,z\right)
&
= \frac{1}{2}\iint e^{-\frac{\lambda}{2}(w-s)^{2}}\chi_{_{0}}(s) s^{2}
e^{-\frac{\lambda}{2}(z-y)^{2}} \theta_{_{0}}(y)
\\
\nonumber
&
\hspace*{11em}
\times
P \left[  \left( \theta_{_{1}}\left(y\right)+  \left(1-\theta_{_{1}}\left(y\right)\right)\right)u(y) \right] \,dsdy
\\
\nonumber
&
=  \frac{1}{2}\iint e^{-\frac{\lambda}{2}(w-s)^{2}}\chi_{_{0}}(s) s^{2}
\left(P^{*}e^{-\frac{\lambda}{2}(z-y)^{2}}\right) \theta_{_{1}}\left(y\right)u(y) \,dsdy
\\
\nonumber
&
\quad
+ \frac{1}{2}\iint e^{-\frac{\lambda}{2}(w-s)^{2}}\chi_{_{0}}(s) s^{2}
\left[P^{*}\left(e^{-\frac{\lambda}{2}(z-y)^{2}} \theta_{_{0}}(y)\right)\right]
\\
\nonumber
&
\hspace*{15em}
\times
\left(1-\theta_{_{1}}\left(y\right)\right) u(y) \,dsdy;
\end{align}
%
	$P^{*}$ denotes the adjoint of $P$.\\

	Case $N=2$; we have
%
\begin{align}\label{N=2}
\mathscr{P}_{2}\left(w,z\right)\doteq
\frac{1}{4!}\iint e^{-\frac{\lambda}{2}(w-s)^{2}} e^{-\frac{\lambda}{2}(z-y)^{2}} \chi_{_{0}}(s)\theta_{_{0}}(y)s^{4}P^{2}u(y) \,dsdy.
\end{align}
%
	We introduce $\theta_{_{1}}(y)$  in $C_{0}^{\infty} \left(\mathscr{B}_{r_{1}}\left(x_{0}\right)\right)$ and $\theta_{_{2}}(y)$ in
	$ C_{0}^{\infty} \left(\mathscr{B}_{r_{2}}\left(x_{0}\right)\right)$ such that
	$\text{supp}\left(\theta_{_{1}}\right)\subseteq \mathscr{B}_{r_{1}}\left(x_{0}\right)$, where $\theta_{_{0}}(y) \equiv 1$, 
	$\theta_{_{1}}(y) \equiv 1$ in $\mathscr{B}_{r_{2}}\left(x_{0}\right)$,
	$\text{supp}\left(\theta_{_{2}}\right)\subseteq \mathscr{B}_{r_{2}}\left(x_{0}\right)$, where $\theta_{_{1}}(y) \equiv 1$,
	and $\theta_{_{2}}(y) \equiv 1$ in $\mathscr{B}_{r_{3}}\left(x_{0}\right)$, where $\tilde{r}_{0} \leq r_{3} < r_{2}< r_{1}<r_{0}$.
	We have
%
\begin{align}\label{N=2_1}
\mathscr{P}_{2}\left(w,z\right)
&
= \frac{1}{4!}\iint e^{-\frac{\lambda}{2}(w-s)^{2}}\chi_{_{0}}(s) s^{4}
e^{-\frac{\lambda}{2}(z-y)^{2}}  \theta_{_{0}}(y)
\\
\nonumber
&
\hspace*{11em}
\times
P \left[ \left( \theta_{_{1}}\left(y\right)+  \left(1-\theta_{_{1}}\left(y\right) \right) \right) Pu(y) \right] \,dsdy
\\
\nonumber
&
=  \frac{1}{4!}\iint e^{-\frac{\lambda}{2}(w-s)^{2}}\chi_{_{0}}(s) s^{4}
\left(P^{*}e^{-\frac{\lambda}{2}(z-y)^{2}}\right) \theta_{_{1}}\left(y\right)Pu(y) \,dsdy
\\
\nonumber
&
\quad
+ \frac{1}{4!}\iint e^{-\frac{\lambda}{2}(w-s)^{2}}\chi_{_{0}}(s) s^{4}
\left[P^{*}\left(e^{-\frac{\lambda}{2}(z-y)^{2}} \theta_{_{0}}(y)\right)\right]
\\
\nonumber
&
\hspace*{17em}
\times
\left(1-\theta_{_{1}}\left(y\right)\right) Pu(y) \,dsdy
\\
\nonumber
&
=  \frac{1}{4!}\iint e^{-\frac{\lambda}{2}(w-s)^{2}}\chi_{_{0}}(s) s^{4}
\left(P^{*}e^{-\frac{\lambda}{2}(z-y)^{2}}\right) \theta_{_{1}}\left(y\right)
\\
\nonumber
&
\hspace*{12em}
\times
P \left[ \left( \theta_{_{2}}\left( y \right) +  \left( 1 - \theta_{_{2}} \left( y \right) \right) \right) u(y) \right] \,dsdy
\\
\nonumber
&
\quad
+ \frac{1}{4!}\iint e^{-\frac{\lambda}{2}(w-s)^{2}}\chi_{_{0}}(s) s^{4}
\left[P^{*}\left(e^{-\frac{\lambda}{2}(z-y)^{2}} \theta_{_{0}}(y)\right)\right]
\\
\nonumber
&
\hspace*{17em}
\times
\left(1-\theta_{_{1}}\left(y\right)\right) Pu(y) \,dsdy
\\
\nonumber
&
=  \frac{1}{4!}\iint e^{-\frac{\lambda}{2}(w-s)^{2}}\chi_{_{0}}(s) s^{4}
\left[\left(P^{*}\right)^{2}e^{-\frac{\lambda}{2}(z-y)^{2}}\right] \theta_{_{2}}\left(y\right) u(y) \,dsdy
\\
\nonumber
&
\quad
+  \frac{1}{4!}\iint e^{-\frac{\lambda}{2}(w-s)^{2}}\chi_{_{0}}(s) s^{4}
\left[ P^{*} \left( \left( P^{*} e^{-\frac{\lambda}{2}(z-y)^{2}}\right) \right.\right.
\\
\nonumber
&
\hspace*{15em}
\left.\left.
\times
\theta_{_{1}} \left( y \right)\right) \right] \left( 1-\theta_{_{2}}\left( y \right) \right) u(y) \,dsdy
\\
\nonumber
&
\quad
+ \frac{1}{4!}\iint e^{-\frac{\lambda}{2}(w-s)^{2}}\chi_{_{0}}(s) s^{4}
\left[P^{*}\left(e^{-\frac{\lambda}{2}(z-y)^{2}} \theta_{_{0}}(y)\right)\right]
\\
\nonumber
&
\hspace*{17em}
\times
\left(1-\theta_{_{1}}\left(y\right)\right) Pu(y) \,dsdy.
\end{align}
%
	The idea is to introduce a sequence of cut-off functions, the support of the subsequent nested where the previous is identically equal to one,
	in order to move all the powers of $ P $ on the exponential function in a neighborhood of $x_{0}$ with the purpose of taking advantage of (\ref{FBI_u}).
	However this will give rise to other terms which still involve powers of $ P $ acting on $u$, but in a region far from $ x_{0}$.  
	We handle the general case as above.  We introduce a family of smooth functions $\lbrace \theta_{_{j}}(y) \rbrace_{1\leq j\leq N}$ such that
	such that $\text{supp}\left(\theta_{_{j}}\right)\subseteq \mathscr{B}_{r_{j}}\left(x_{0}\right)$
	and $\theta_{_{j}}(y) \equiv 1$ in $\mathscr{B}_{r_{j+1}}\left(x_{0}\right)$, $\tilde{r}_{0} \leq r_{N+1}< r_{N}< \,\dots\, < r_{2}< r_{1}<r_{0}$.
	So, we see that for every $j$ less or equal than $N+1$, one has that: $ i $ less than $ j $ implies $\theta_{_{i}}\equiv 1$ on a neighborhood of $ \theta_{_{j}}$. 
	In order to construct the functions $\theta_{j}$ we follow the same strategy used to construct
	the  Ehrenpreis-H\"or\-man\-der cut-off functions. More precisely we choose $r_{j} =r_{0}- \left(r_{0}-\tilde{r}_{0}\right) \frac{j}{N+1}$,
	we have $r_{j} - r_{j+1} = \frac{r_{0}-\tilde{r}_{0}}{N+1}$. Let $\psi$ be a function in $ \mathscr{D}(\mathbb{R}^{n})$ with support in
	$\mathscr{B}_{1/4}(0) \doteq \lbrace y\in \mathbb{R}^{n}\,: \, |y| \leq 1/4 \rbrace$ such that $\psi \geq 0$ and
	$\int \psi \, dy =1$. For every $\gamma > 0$ we write $\psi_{\gamma}(y) = \gamma^{-n} \psi\left(\frac{x}{\gamma}\right)$.
	Let $\chi_{j}$ be the characteristic function of the set
	$\lbrace y\in \mathbb{R}^{n}\,:\, \text{dist}\left(y; \mathscr{B}_{r_{j+1}}(x_{0})\right) < \frac{r_{0}-\tilde{r}_{0}}{2(N+1)}\rbrace$.
	We set
%
\begin{align*}
\theta_{j} = \psi_{\frac{r_{0}-\tilde{r}_{0}}{\left(N+1\right)}} * \psi_{\frac{r_{0}-\tilde{r}_{0}}{\left(N+1\right)}} *\,  \chi_{j}.
\end{align*}
%
	These functions have the desired properties. Moreover we have
%
\begin{align*}
&\| D_{y_{i}} \theta_{j} \|_{\infty}
\leq \|D_{y_{i}}\psi_{\frac{r_{0}-\tilde{r}_{0}}{\left(N+1\right)}}\|_{L^{1}} \|\psi_{\frac{r_{0}-\tilde{r}_{0}}{\left(N+1\right)}}\|_{L^{1}}\|\chi_{j}\|_{\infty}
\leq C_{0}\frac{N+1}{r_{0}-\tilde{r}_{0}},
\\
&
\| D_{y_{i}}D_{y_{k}}  \theta_{j}\|_{\infty} \leq 
\|D_{y_{i}}\psi_{\frac{r_{0} -\tilde{r}_{0}}{\left(N+1\right)}}\|_{L^{1}} \|D_{y_{k}}\psi_{\frac{r_{0} -\tilde{r}_{0}}{\left(N+1\right)}}\|_{L^{1}} \|\chi_{j}\|_{\infty}
\leq \left( C_{0}\frac{N+1}{r_{0}-\tilde{r}_{0}}\right)^{2},
\end{align*}
%
	where $C_{0} = \sup_{1\leq i\leq n} \|D_{y_{i}} \psi  \|_{L^{1}\left(\mathscr{B}_{1/4}(0)\right)}$.
%
\begin{remark}
	One may also choose a sequence of $\theta_{j}$ independent of $N$, by repeating the above construction 
	and taking the convolution with   $\psi_{( r_{0} - \tilde{r}_{0} )/2^{j} }$.\\
	Moreover the $\theta_{_{j}}$ can be constructed by just one convolution,
	i.e. $\theta_{_{j}} = \psi_{\frac{r_{0}-\tilde{r}_{0}}{ 2^{j} }} *\,  \chi $.
	This will be more evident in the next few steps.
\end{remark}
\noindent
	We set $X_{j} (x,D) = \displaystyle\sum_{\ell=1}^{n} a_{\ell,j}(x) D_{i}$.
	We have
%
\begin{align}\label{case_N}
&\mathscr{P}_{N}(w,z)
=\frac{1}{(2N)!}\iint e^{-\frac{\lambda}{2}(w-s)^{2}} e^{-\frac{\lambda}{2}(z-y)^{2}} \chi_{_{0}}(s)\theta_{_{0}}(y)s^{2N}P^{N}u(y) \,ds dy
\\
\nonumber
&
=\frac{1}{(2N)!}\iint e^{-\frac{\lambda}{2}(w-s)^{2}}  \chi_{0}(s)s^{2N}
\left[ \left(P^{*}\right)^{N}e^{-\frac{\lambda}{2}(z-y)^{2}}\right] \theta_{_{N}}(y)u(y) \,ds dy
\\
\nonumber
&\,\,\,+\frac{1}{(2N)!} \sum_{j=1}^{N}
\iint e^{-\frac{\lambda}{2}(w-s)^{2}}\! \chi_{0}(s)s^{2N}\!\!
\left\{ \! P^{*} \! \left[ \! \left( \left( P^{*}\right)^{j-1}e^{-\frac{\lambda}{2}(z-y)^{2}} \right) \theta_{j-1}(y) \right]\right\}
\\
\nonumber
&
\hspace*{19em}
\times\left( 1- \theta_{j}(y)\right) P^{N-j}u(y) \,ds dy 
\\
\nonumber
&
=\frac{1}{(2N)!}\iint e^{-\frac{\lambda}{2}(w-s)^{2}}  \chi_{0}(s)s^{2N}
\left[ \left(P^{*}\right)^{N}e^{-\frac{\lambda}{2}(z-y)^{2}}\right] \theta_{_{N}}(y)u(y) \,ds dy
\\
\nonumber
&\,\,\,+\frac{1}{(2N)!} \sum_{j=1}^{N}
\iint e^{-\frac{\lambda}{2}(w-s)^{2}}\! \chi_{0}(s)s^{2N}
\! \left( \left( P^{*}\right)^{j}e^{-\frac{\lambda}{2}(z-y)^{2}} \right) \theta_{j-1}(y) 
\\
\nonumber
&
\hspace*{19em}
\times
\left( 1- \theta_{j}(y)\right) P^{N-j}u(y) \,ds dy 
\\
\nonumber
&\,\,\,+\frac{1}{(2N)!} \sum_{j=1}^{N}
\iint e^{-\frac{\lambda}{2}(w-s)^{2}}\! \chi_{0}(s)s^{2N}\!
\left( \left( P^{*}\right)^{j-1}e^{-\frac{\lambda}{2}(z-y)^{2}} \right) \left( P\theta_{j-1}(y) \right)
\\
\nonumber
&
\hspace*{19em}
\times
\left( 1- \theta_{j}(y)\right) P^{N-j}u(y) \,ds dy 
\\
\nonumber
&\,\,\,+\frac{1}{(2N)!} \sum_{j=1}^{N}
\iint e^{-\frac{\lambda}{2}(w-s)^{2}}\! \chi_{0}(s)s^{2N}
\left[ \,\sum_{k=1}^{m}
\left( X_{k}\left( P^{*}\right)^{j-1}e^{-\frac{\lambda}{2}(z-y)^{2}} \right)\right.
\\
\nonumber
&
\hspace*{13em}
\times
\left( X_{k}\theta_{j-1}(y) \right)\bigg]
\left( 1- \theta_{j}(y)\right) P^{N-j}u(y) \,ds dy 
%
\\
\nonumber
&\,\,\,+\frac{2}{(2N)!} \sum_{j=1}^{N}
\iint e^{-\frac{\lambda}{2}(w-s)^{2}}\! \chi_{0}(s)s^{2N} 
\left[ \left(\left( P^{*}\right)^{j-1}e^{-\frac{\lambda}{2}(z-y)^{2}} \right)\right.
\\
\nonumber
&
\hspace*{11em}
\times
\left.
\sum_{k=1}^{m} f_{k}
\left( X_{k}\theta_{j-1}(y) \right)\right]
\!\!\left( 1- \theta_{j}(y)\right) P^{N-j}u(y) \,ds dy 
\\
\nonumber
& =  I_{1} + I_{2} + I_{3} + I_{4} + I_{5},
\end{align}
%
	where 
	$f_{k} = \frac{1}{i} \sum_{\ell= 1}^{n} a_{\ell,k}^{(e_{\ell})}(y)$,
	$e_{\ell}$, $\ell, =1,\dots,n$, in the upper index denotes the derivatives in the direction $\ell$.\\ 
	Before estimating the above terms a few remarks are in order:
%
\begin{enumerate}
	%
		\item[i)]  each step no more than two derivatives act on $\theta_{_{j}}(y)$;
		\item[ii)] let $ y \in \text{supp}(\theta_{_{j-1}}^{(\alpha)})\cap \text{supp}\left(1-\theta_{_{j}}\right)$, $ 0 \leq |\alpha| \leq 2 $, since
		$\Omega $ is a complex neighborhood of $x_{0}-i\xi_{0}$ such that $\pi_{z'} \left( \Omega\right) \subset \mathscr{B}_{\tilde{r}_{0} -\varepsilon}$, 
		we have that $(z'-y)^{2}\geq \varepsilon^{2}$;
		\item[iii)] without loss of generality we may write
	%
	\begin{align*}
	\left(P^{*}(y,D)\right)^{N} = \sum_{ |\beta| \leq 2N} a_{_{2N,\beta}}(y) D^{\beta} 
	\end{align*}  
	%
		where $a_{_{2N,\beta}}(y)$ are analytic functions such that for any compact set $K$ in $U$ we have
	%
	\begin{align}\label{est_a_2N_beta}
	\left| a_{_{2N,\beta}}^{(\gamma)} (y)\right|  \leq C_{K}^{3N - |\beta|+|\gamma|} \left(2N- |\beta| + |\gamma|\right)!
	\qquad \forall \, y \in K \text{ and } \gamma \in \mathbb{Z}^{n}_{+}.
	\end{align}
	%
		\item[iv)] the following identity holds
	%
	\begin{align*}
	&\left(\frac{d}{d y_{k}}\right)^{\beta_{k}} \!\!\! e^{-\frac{\lambda}{2}\left( z_{k}-y_{k}\right)^{2}} 
	=
	e^{-\frac{\lambda}{2}\left( z_{k}-y_{k}\right)^{2}}  \sum_{\ell_{k} =0}^{\lfloor \frac{\beta_{k}}{2}\rfloor}
	\frac{\beta_{k}!(i)^{2(\beta_{k}-\ell_{k})}}{\ell_{k}! \left( \beta_{k} -2 \ell_{k}\right)! 2^{\ell_{k}}} 
	\lambda^{\beta_{k}-\ell_{k}} \left( z_{k} -y_{k}\right)^{\beta_{k}-2 \ell_{k}}
	\\
	&
	\quad=e^{-\frac{\lambda}{2}\left( z_{k}-y_{k}\right)^{2}} (i)^{\beta_{k}} \left(\frac{\lambda}{2}\right)^{\beta_{k}/2}
	\sum_{\ell_{k} =0}^{\lfloor \frac{\beta_{k}}{2}\rfloor}
	\frac{\beta_{k}!}{\ell_{k}! \left( \beta_{k} -2 \ell_{k}\right)!} \left[i\sqrt{2\lambda}\left( z_{k} -y_{k}\right) \right]^{\beta_{k}-2 \ell_{k}}.
	\end{align*}
\end{enumerate}
%
\vspace{3em}
	\textbf{Estimate of the term $I_{2}$}. Since we are far from $x_{0}$ we expect exponential decay.
	We have
%
\begin{align*}
I_{2}
&= \frac{1}{(2N)!} \sum_{j=1}^{N}
\int e^{-\frac{\lambda}{2}(w-s)^{2}}\! \chi_{0}(s)s^{2N} ds  \int\sum_{ |\beta| \leq 2j} \frac{1}{(i)^{|\beta|}} a_{_{2j,\beta}} (y)
e^{\frac{\lambda}{2}(z'')^{2}+ i \lambda(y-z')z''}
\\
& \times
\left[\prod_{\nu=1}^{n}  \left( \sum_{\gamma_{\nu} =0}^{\lfloor\frac{\beta_{\nu}}{2}\rfloor}  
\frac{\beta_{\nu}! i^{2(\beta_{\nu}-\gamma_{\nu})}}{\gamma_{\nu}! \left( \beta_{\nu}-2\gamma_{\nu}\right)! 2^{|\gamma_{\nu}|}} 
\left(\lambda^{\beta_{\nu} - \gamma_{\nu}} \left( z_{\nu}- y_{\nu}\right)^{\beta_{\nu} - 2\gamma_{\nu}}\right)
e^{-\frac{\lambda}{2}(z_{\nu}-y_{\nu})^{2}}
\right)
\right]
\\
\vspace*{0.5em}
& \times
\theta_{j-1}(y)\left( 1- \theta_{j}(y)\right) P^{N-j}u(y) \, dy. 
\end{align*}
%
	We remark that the integral with respect the variable $ s $ is the FBI transform of $ \chi_{0}(s)s^{2N}$.
	We take $\Re(w) \in ]-\delta_{2}-\sqrt{\tilde{\varepsilon}_{0}}, \delta_{2} + \sqrt{\tilde{\varepsilon}_{0}}[ $, $\tilde{\varepsilon}_{0} $ sufficiently small positive constant.
	Splitting the domain of integration in the regions where $\chi_{0}(s)\neq 1$  and $\chi_{0}(s)=1$ and changing, in the last one region, the integration path
	as in the Remark \ref{RK_def_2}, so that it is in the strip $\sigma = s+ i \sigma''$, $|\sigma''| <\delta_{2}/2$, where we consider the holomorphic extension of $s^{2N}$,
	we can conclude that there is a positive constants $A $ such that
%
\begin{align*}
\left|
\int e^{-\frac{\lambda}{2}(w-s)^{2}}\! \chi_{0}(s)s^{2N} ds 
\right|
\leq e^{\frac{\lambda}{2} \left( w''\right)^{2}} A \delta_{1}^{ 2N} e^{-\lambda\tilde{\varepsilon}_{0}}.
\end{align*}
%
	Since $ y \in \mathscr{B}_{r_{j-1}}(x_{0})\setminus \mathscr{B}_{r_{j+1}}\left(x_{0}\right)$ we have $(z'-y)^{2} \geq \varepsilon_{0}$. 
	We obtain
%
\begin{align*}
&\left|\lambda^{\beta_{\nu} - \gamma_{\nu}} \left( z_{\nu}- y_{\nu}\right)^{\beta_{\nu} - 2\gamma_{\nu}}e^{-\frac{\lambda}{2}(z'_{\nu}-y_{\nu})^{2}} \right|
\\
&\leq 
2\cdot 2^{\frac{3}{2}(\beta_{\nu} - 2\gamma_{\nu})} \left(\frac{8}{\varepsilon_{0}}\right)^{\beta_{\nu} - \gamma_{\nu}}
\left(\beta_{\nu} !\right)^{\frac{1}{2}} \left[\left(\beta_{\nu} - 2\gamma_{\nu}\right)!\right]^{\frac{1}{2}}
e^{-\frac{\varepsilon_{0}}{16}\lambda},
\end{align*}
%
	where we can assume that $ | z''_{\nu}| \leq 1$.
	Since $\beta_{\nu} ! \leq 2^{\beta_{\nu} + 2\gamma_{\nu}} \left[ \left(\beta_{\nu} - 2\gamma_{\nu}\right)!\right] \left( \gamma_{\nu}!\right)^{2}$,
	we have
%
\begin{align*}
\sum_{\gamma_{\nu} =0}^{\lfloor\frac{\beta_{\nu}}{2}\rfloor}  
\frac{\beta_{\nu}! }{\gamma_{\nu}! \left( \beta_{\nu}-2\gamma_{\nu}\right)! 2^{|\gamma_{\nu}|}} 
\left|\left(\lambda^{\beta_{\nu} - \gamma_{\nu}} \left( z_{\nu}- y_{\nu}\right)^{\beta_{\nu} - 2\gamma_{\nu}}\right)
e^{-\frac{\lambda}{2}(z_{\nu}-y_{\nu})^{2}}\right|
\\
\leq 4\cdot \left( \beta_{\nu}! \right) \left(\frac{32}{\varepsilon_{0}}\right)^{\beta_{\nu} }
e^{-\frac{n\varepsilon_{0}}{16}\lambda},
\end{align*}
%
	then
%
\begin{align*}
\prod_{\nu=1}^{n}  \left( \sum_{\gamma_{\nu} =0}^{\lfloor\frac{\beta_{\nu}}{2}\rfloor}  
\frac{\beta_{\nu}!}{\gamma_{\nu}! \left( \beta_{\nu}-2\gamma_{\nu}\right)! 2^{|\gamma_{\nu}|}} 
\left|\left(\lambda^{\beta_{\nu} - \gamma_{\nu}} \left( z_{\nu}- y_{\nu}\right)^{\beta_{\nu} - 2\gamma_{\nu}}\right)
e^{-\frac{\lambda}{2}(z_{\nu}-y_{\nu})^{2}}\right|
\right)
\\
\leq 4^{n}  \left( \beta !\right) \left(\frac{32}{\varepsilon_{0}}\right)^{|\beta| }
e^{-\frac{n\varepsilon_{0}}{16}\lambda}.
\end{align*}
%
	We obtain
%
\begin{align*}
&| e^{-\lambda \phi_{0}(w,z)} I_{2}| 
\leq 4^{n} A\, \frac{(2\pi)^{n}r_{0}^{n-1}}{\Gamma\left(\frac{n}{2}\right)}\,
e^{-\tilde{\varepsilon}_{0}\lambda} \, e^{-\frac{n\varepsilon_{0}}{16}\lambda}\delta_{1}^{N}
\\
&
\times
\sum_{j=1}^{N} \frac{1}{\left(2N\right)!}
\left(  \sum_{ |\beta| \leq 2j} C_{1}^{3j-|\beta|+1} \left(2j-|\beta|\right)! \beta!  \left(\frac{32}{\varepsilon_{0}}\right)^{|\beta| }\right)
\tilde{C}_{2}^{2(N-j)+1} \left[2(N-j)\right]!,
\end{align*}
%
	where $C_{1}$ and $\tilde{C}_{2}$ are the constants in (\ref{est_a_2N_beta}) and in (\ref{An_v_est}), respectively,
	with $K= \overline{\mathscr{B}_{r_{0}}(x_{0})}$. Without loss of generality we may assume that $C_{1}$ and $\tilde{C}_{2}$
	are greater then $2$. Since $ \left( 2j -|\beta|\right)! \leq (2j)! \left( |\beta|!\right)^{-1}$ and $\left[2(N-j)\right]! \leq (2N)! \left( (2j)!\right)^{-1}$,
	we have
%
\begin{align*}
| e^{-\lambda \phi_{0}(w,z)} I_{2}| 
\leq 2 \cdot 8^{n}\, \frac{(2\pi)^{n}r_{0}^{n-1}}{\Gamma\left(\frac{n}{2}\right)}\,
A C_{1} \tilde{C}_{2}\,
e^{-\tilde{\varepsilon}_{0}\lambda} \, e^{-\frac{n\varepsilon_{0}}{16}\lambda}
\delta_{1}^{N} \left( \frac{32\, C_{1}^{\frac{3}{2}} \tilde{C}_{2}}{\varepsilon_{0}} \right)^{2N}.
\end{align*}
%
	Taking $\delta_{1}$ small enough we conclude that there are two positive constants $ C_{2}$  and $\varepsilon_{2}$, independent by $N$, such that
%
\begin{align}\label{Est-I_2}
| e^{-\lambda \phi_{0}(w,z)} I_{2}| \leq C_{2} \left( \frac{1}{2}\right)^{N} e^{-\varepsilon_{2}\lambda}.
\end{align}
%
\vspace{3em}

\noindent
	\textbf{Estimate of the terms $I_{3}$, $I_{4}$ and $I_{5}$}.
	The only difference from $I_{2}$ is that either two derivatives or one derivative act on the functions $\theta_{j}(y)$.
	These terms are treated analogously to the term $I_{2}$.
	Then there are positive constants, $C_{3}$, $C_{4}$ and $C_{5}$, independent of $N$, such that
%
\begin{align}\label{Est-I_3}
| e^{-\lambda \phi_{0}(w,z)} I_{3}| \leq C_{3} (N+1)^{2}  \left( \frac{1}{2}\right)^{N} e^{-\varepsilon_{2}\lambda},
\end{align}
%
	and
%
\begin{align}\label{Est-I_4}
| e^{-\lambda \phi_{0}(w,z)} I_{4}| \leq C_{4} (N+1) \left( \frac{1}{2}\right)^{N} e^{-\varepsilon_{2}\lambda}
\end{align}
%
	and
%
\begin{align}\label{Est-I_5}
| e^{-\lambda \phi_{0}(w,z)} I_{4}| \leq C_{5} (N+1) \left( \frac{1}{2}\right)^{N} e^{-\varepsilon_{2}\lambda}.
\end{align}

\vspace{0.5em}
\noindent
	\textbf{Estimate of the term $I_{1}$}. Roughly speaking we are studying the micro-local regularity of the product of an analytic function
	with $u$ at the point $(x_{0},\xi_{0})$.  In order to estimate this term we take advantage from the following theorem 
	which characterizes micro-local smoothness in terms of $(s_{0}-1)$-almost analytic extendability in certain wedges.
%
\begin{theorem}[see  Theorem 2.3 in \cite{Ber_Hai_2017}]
	Let $u \in \mathscr{D}'\left(U\right)$. Then $\left( x_{0}, \xi_{0}\right) \notin WF_{s_{0}}(u)$ if and only if there exist a neighborhood
	$U_{0}$ of $x_{0}$, open acute cones $\Gamma^{1}, \dots ,\, \Gamma^{k}$ in $\mathbb{R}^{n}\setminus \lbrace 0 \rbrace$
	and $(s_{0}-1)$-almost analytic functions $f_{j}$ on $U_{0} + i \Gamma^{j}_{\varepsilon_{1}}$,
	$\Gamma_{\varepsilon_{1}}^{j}= \Gamma^{j}\cap\lbrace \xi: \, |\xi| <\varepsilon_{1} \rbrace$,
	of temperate growth such that $u =\sum_{j=1}^{k} bf_{j}$ near $x_{0}$ and $\xi_{0} \cdot \Gamma^{j} < 0$ for all $j$.
\end{theorem}
%
	Analogous results in the smooth and analytic category can be found in \cite{BCH_book} and \cite{BH_2012}.
	We point out that in the analytic case the $f_{j}$ are holomorphic functions.
	We recall
%
\begin{definition}\label{Def-almost}
	Let $f \in G^{s_{0}}\left( U\right)$, $U$ open subset of $\mathbb{R}^{n}$, and suppose $\tilde{U}$ is a open neighborhood
	of $U$ in $\mathbb{C}^{n}$. A function $\tilde{f}(y,\eta) \in C^{\infty} ( \tilde{U})$ is called an $(s_{0}-1)$-almost analytic extension
	of $f$ if $\tilde{f}(y,0) = f(y)\, \forall \, y \in U$ and for every compact $K$ in $U$ there exists positive constants $C_{K}$ and small $\varepsilon_{_{K}}$
	such that 
	\begin{align*}
	\left| \partial_{\bar{z}_{j}}\tilde{f} \right| \leq C_{K} e^{-\varepsilon_{_{K}} |\eta|^{-\frac{1}{s_{0}-1}}}, \qquad j=1,\dots,n, 
	\end{align*}  
	holds for $ y \in K $ and $\eta$ in ball of radius $ \varepsilon_{_{K}} $.
\end{definition}
%
	The $(s_{0}-1)$-almost analytic extension of a Gevrey function $f $ can be obtained in the following way
%
\begin{align*}
\tilde{f}(y+i\eta) = \sum_{ \gamma} f^{(\gamma)}(y)\frac{i^{|\gamma|}y^{\gamma}}{\gamma!} \Theta\left( \tilde{c} |\gamma|^{s_{0}-1} |\eta|\right),
\end{align*}
%
	where $\Theta$ is in $C^{\infty}_{0}\left(\mathbb{R}\right)$ such that $\text{supp} \Theta \subset \left[ -1,1\right]$ and $\Theta(y) \equiv 1$
	on $\left[ -1/2, 1/2\right]$. For other details see \cite{Ch_PHD_t} or \cite{Ber_Hai_2017}. We point out that, by hypothesis,
	we can construct in a suitable region an $(s_{0}-1)$-almost analytic extension of $u$.
	We have to estimate
%
\begin{align}\label{Term_I-1}
I_{1}= 
\frac{1}{(2N)!} 
&\int e^{-\frac{\lambda}{2}(w-s)^{2}}  \chi_{0}(s)s^{2N} ds \sum_{ |\beta| \leq 2N}\sum_{\substack{\gamma \leq \lfloor \frac{\beta}{2}\rfloor \\  \gamma \in \mathbb{Z}^{n}_{+} }}
\frac{\beta!\, i^{|\beta|-2|\gamma|}}{\gamma! \left( \beta-2\gamma!\right) 2^{\gamma}}
\lambda^{|\beta|-|\gamma|} 
\\
\nonumber
&
\times
\int a_{_{2N,\beta}} (y)
e^{-\frac{\lambda}{2}(z-y)^{2}}
\prod_{\nu=1}^{n} \left( z_{\nu}- y_{\nu}\right)^{\beta_{\nu} - 2\gamma_{\nu}}
\theta_{_{N}}(y)u(y) \, dy.
\end{align}
%
	In order to handle the integral with respect the variable $y$ we follow the classical strategy developed by 
	Bros and Iagolnitzer, \cite{Iag_75}. We split the integration domain in two parts:
	$\mathscr{B}_{r_{N}}(x_{0}) \setminus \mathscr{B}_{\tilde{r}_{0}}(x_{0})$ and $ \mathscr{B}_{\tilde{r}_{0}}(x_{0})$.
	We have
%
\begin{align}\label{Split_int}
& \int a_{_{2N,\beta}} (y)
e^{-\frac{\lambda}{2}(z-y)^{2}}
\prod_{\nu=1}^{n} \left( z_{\nu}- y_{\nu}\right)^{\beta_{\nu} - 2\gamma_{\nu}}
\theta_{_{N}}(y)u(y) \, dy
\\
&
\nonumber
\quad= \int_{\mathscr{B}_{r_{N}}(x_{0}) \setminus \mathscr{B}_{\tilde{r}_{0}}(x_{0})} a_{_{2N,\beta}} (y)
e^{-\frac{\lambda}{2}(z-y)^{2}}
\prod_{\nu=1}^{n} \left( z_{\nu}- y_{\nu}\right)^{\beta_{\nu} - 2\gamma_{\nu}}
\theta_{_{N}}(y)u(y) \, dy
\\
&
\nonumber
\qquad
+
\int_{\mathscr{B}_{\tilde{r}_{0}}(x_{0})} a_{_{2N,\beta}} (y)
e^{-\frac{\lambda}{2}(z-y)^{2}}
\prod_{\nu=1}^{n} \left( z_{\nu}- y_{\nu}\right)^{\beta_{\nu} - 2\gamma_{\nu}}u(y) \, dy.
\end{align}
%
	In the first region $\left( z' - y \right)^{2} \geq \varepsilon_{0}$, this will give the analytic exponential decay in this region. 
	Our purpose is to verify that the second integral gives the desired Gevrey exponential decay.
	Since $\left( x_{0}, \xi_{0}\right) \notin WF_{s_{0}}(u)$ without loss of generality we may assume that
	$u$ is a boundary value of $\tilde{u}(\zeta)$, $(s_{0}-1)$-almost analytic function on $\mathscr{B}_{\tilde{r}_{0}}(x_{0}) + i \Gamma_{\varepsilon_{2}} $,
	$\Gamma_{\varepsilon_{2}} = \lbrace \eta \in \Gamma:\, |\eta| < \varepsilon_{2} \rbrace$, where $\Gamma$ is an open cone such that
	$\eta\cdot \xi_{0} < 0$ for all $\eta \in \Gamma$. We point out that, for a fixed a neighborhood of $\xi_{0}$, we can choose $\Gamma$ such that
	$\eta\cdot \xi < 0$ for all $\eta \in \Gamma$ and $\xi$ in the neighborhood of $\xi_{0}$.
	On the other hand $a_{_{2N,\beta}} (y)$ are analytic functions, we can construct their holomorphic extension $\tilde{a}_{_{2N,\beta}} (\zeta)$
	in $\mathbb{C}^{n}_{\zeta}$, where $\zeta = y+i \eta$ and $|\eta | \leq \varepsilon_{3}$.
	We take $\varepsilon_{2}$ such that $\varepsilon_{2} \leq \varepsilon_{3}$.\\
	Let $\vartheta(y) \in C^{\infty}_{0} \left( \mathbb{R}^{n}\right)$ with support equal to $\mathscr{B}_{\tilde{r}_{0}}\left( x_{0} \right)$
	such that $0 \leq \vartheta(y) \leq 1$ and $\vartheta(x_{0}) = 1$. Let $\eta^{0} \in \Gamma_{\varepsilon_{2}}$, we define
	the $n$-dimensional manifold, $S_{\eta^{0},\varepsilon_{4}}$, in $\mathbb{C}^{n}_{\zeta}$, $\zeta = y +i \eta$, given by 
%
\begin{align*}
y \mapsto \zeta  = y +i\varepsilon_{4} \vartheta(y)\eta^{0},
\end{align*} 
%
	where $\varepsilon_{4} \in \mathbb{R}_{+}$ and is sufficiently small so that  $S_{\eta^{0},\varepsilon_{4}}$
	is contained in $\mathscr{B}_{\tilde{r}_{0}}(x_{0}) + i \Gamma_{\varepsilon_{2}} $. We remark that the boundary
	of $S_{\eta^{0},\varepsilon_{4}}$ is equal to $\partial\mathscr{B}_{\tilde{r}_{0}}(x_{0})$.
	By the Stokes theorem we have
%
\begin{align*}
& \int_{\mathscr{B}_{\tilde{r}_{0}}(x_{0})} a_{_{2N,\beta}} (y)
e^{-\frac{\lambda}{2}(z-y)^{2}}
\prod_{\nu=1}^{n} \left( z_{\nu}- y_{\nu}\right)^{\beta_{\nu} - 2\gamma_{\nu}}u(y) \, dy
\\
&
\qquad
= - \int_{S_{\eta^{0},\varepsilon_{4}}} \tilde{a}_{_{2N,\beta}} (\zeta)
e^{-\frac{\lambda}{2}(z-\zeta)^{2}}
\prod_{\nu=1}^{n} \left( z_{\nu}- \zeta_{\nu}\right)^{\beta_{\nu} - 2\gamma_{\nu}}\tilde{u}(\zeta) \, d\zeta
\\
&
\qquad\qquad
+ \int_{D_{\eta^{0}}} 
d\left( \tilde{a}_{_{2N,\beta}} (\zeta)
e^{-\frac{\lambda}{2}(z-\zeta)^{2}}
\prod_{\nu=1}^{n} \left( z_{\nu}- \zeta_{\nu}\right)^{\beta_{\nu} - 2\gamma_{\nu}}\tilde{u}(\zeta) \right)\, \wedge d\zeta,
\end{align*} 
%
	where $D_{\eta^{0}} = \displaystyle\bigcup_{0 < t < \varepsilon_{4}} S_{\eta^{0},t } \subseteq \mathscr{B}_{\tilde{r}_{0}}(x_{0}) + i \Gamma_{\varepsilon_{2}}$
	and $\partial D_{\eta^{0}} = \mathscr{B}_{\tilde{r}_{0}}(x_{0}) \cup S_{\eta^{0},\varepsilon_{4} }$.\\
	Since
%
\begin{align*}
&d\zeta_{j} = \sum_{\substack{k=1\\  k\neq j}}^{n}(i t\, \vartheta^{(e_{k})} (y)\,\eta^{0}_{j}) dy_{k} 
+ (1 + i t\, \vartheta^{(e_{j})}(y)\,\eta^{0}_{j}) dy_{j}+ ( i \,\vartheta(y)\,\eta_{j}^{0})dt,
\\
&
d\bar{\zeta_{i}}= \sum_{\substack{k=1\\  k\neq j}}^{n}(-i t\, \vartheta^{(e_{k})} (y)\,\eta^{0}_{j}) dy_{k} 
+ (1 - i t\, \vartheta^{(e_{j})}(y)\,\eta^{0}_{j}) dy_{j} - ( i \,\vartheta(y)\,\eta_{j}^{0})dt,
\end{align*}
%
	where $ \vartheta^{(e_{j})}(y) = \left( \partial_{y_{j}}\vartheta\right)(y)$,
	and $\tilde{a}_{_{2N,\beta}}(\zeta)$ are holomorphic functions, analytic extensions of $a_{_{2N,\beta}}(y)$, in $D_{\eta^{0}}$,
	we have 
%
\begin{align*}
&\left(d\left( \tilde{a}_{_{2N,\beta}} (\zeta)
e^{-\frac{\lambda}{2}(z-\zeta)^{2}}
\prod_{\nu=1}^{n} \left( z_{\nu}- \zeta_{\nu}\right)^{\beta_{\nu} - 2\gamma_{\nu}}\tilde{u}(\zeta) \right)\, \wedge d\zeta\right)_{|_{S_{\eta^{0},t }}}
\\
&=
\sum_{j=1}^{n}\left( \tilde{a}_{_{2N,\beta}} (\zeta ) e^{-\frac{\lambda}{2}(z-\zeta)^{2}}
\prod_{\nu=1}^{n} \left( z_{\nu}- \zeta_{\nu} \right)^{\beta_{\nu} - 2\gamma_{\nu}} 
\frac{\partial \tilde{u}}{\partial \bar{\zeta}_{j} }(\zeta)  d \bar{\zeta_{j}} \wedge  d\zeta \right)_{|_{S_{\eta^{0},t }}}
\\
&
=
\sum_{j=1}^{n}\left(\! \tilde{a}_{_{2N,\beta}} (\zeta ) e^{-\frac{\lambda}{2}(z-\zeta)^{2}}
\prod_{\nu=1}^{n} \left( z_{\nu}- \zeta_{\nu} \right)^{\beta_{\nu} - 2\gamma_{\nu}} 
\frac{\partial \tilde{u}}{\partial \bar{\zeta}_{i} }(\zeta) \right)_{|_{S_{\eta^{0},t }}}
\hspace*{-2em}
\det\left(A_{j}(y,t,\eta^{0})\right)\,dt\,dy,
\end{align*}
%
	where $A_{j}(y,t,\eta^{0}) $ is the $(n+1)\times(n+1)-$matrix
%
\vspace{0.3em}
\begin{equation*}
\hspace*{-4em}
\tiny{
	\begin{pmatrix}
	-i t \vartheta^{(e_{1})} (y)\eta^{0}_{j} \!\!\! &\!\!\! \cdots\!\!\!  &\!\!\!  -i t \vartheta^{(e_{j-1})} (y)\eta^{0}_{j}&\!\!\!  1 -i t \vartheta^{(e_{j})} (y)\eta^{0}_{j}&
	\!\!\! -i t \vartheta^{(e_{j+1})} (y)\eta^{0}_{j} & \!\!\! \cdots &\!\!\!  -i  \vartheta (y)\eta^{0}_{j} 
	\\
	1 +i t\, \vartheta^{(e_{1})} (y)\eta^{0}_{1} \!\!\! & \!\!\! i t \vartheta^{(e_{2})} (y)\eta^{0}_{1}\!\!  & \!\!\! \cdots & \!\!\!  \cdots & \!\!\!  \cdots
	&\!\!\!   \cdots &\!\!\!  i  \vartheta (y)\eta^{0}_{1} 
	\\
	i t \vartheta^{(e_{1})} (y)\eta^{0}_{1} \!\!\! & \!\!\! 1 +i t \vartheta^{(e_{2})} (y)\eta^{0}_{2} & \!\! i t \vartheta^{(e_{3})} (y)\eta^{0}_{2}\!\! 
	&\!\!\!   \cdots \!\!\! &\!\!\!   \cdots \!\!\! &\!\!\!   \cdots \!\!\! &\!\!\!  i  \vartheta (y)\eta^{0}_{2} \\
	\vdots  & \!\!\!  \ddots & \!\!\! \ddots & \!\!\! \ddots &\!\!\!  \ddots &\!\!\!  \ddots &\!\!\!  \vdots  \\
	i t \vartheta^{(e_{1})} (y)\eta^{0}_{n} \!\!\! & \!\!\!  \cdots\!\!\!  & \!\!\! \cdots\!\!\!  & \!\!\! \cdots\!\!\!  &\!\!\!  i t \vartheta^{(e_{n-1})} (y)\eta^{0}_{n}
	&\!\!\!  1+  i t \vartheta^{(e_{n})} (y)\eta^{0}_{n}
	&\!\!\!   i  \vartheta (y)\eta^{0}_{n} 
	\end{pmatrix}
}
\end{equation*}

\vspace{1em}

	We obtain 
%
\begin{align*}
e^{-\lambda \phi_{0}(w,z)} I_{1}=  e^{-\lambda \phi_{0}(w,z)} I_{1,1} + e^{-\lambda \phi_{0}(w,z)} I_{1,2} + e^{-\lambda \phi_{0}(w,z)} I_{1,3},
\end{align*}
%
	where
%
\begin{align*}
&e^{-\lambda \phi_{0}(w,z)}  I_{1,1}=
\frac{1}{(2N)!}  e^{-\frac{\lambda}{2} \left[\left( w''\right)^{2} + \left( z''\right)^{2}\right]} 
\int e^{-\frac{\lambda}{2}(w-s)^{2}}  \chi_{0}(s)s^{2N} ds
\\
&
\times
\quad
\sum_{ |\beta| \leq 2N}\sum_{\substack{\gamma \leq \lfloor \frac{\beta}{2}\rfloor \\  \gamma \in \mathbb{Z}^{n}_{+} }}
\frac{\beta!\, i^{|\beta|-2|\gamma|}}{\gamma! \left( \beta-2\gamma!\right) 2^{\gamma}}
\lambda^{|\beta|-|\gamma|} 
\!\int_{\mathscr{B}_{r_{N}}(x_{0}) \setminus \mathscr{B}_{\tilde{r}_{0}}(x_{0})}
\!\!\!\!\!\!\!\!\!\!\! e^{-\frac{\lambda}{2}(z-y)^{2}} a_{_{2N,\beta}} (y) 
\\
&
\hspace{15em}
\times
\prod_{\nu=1}^{n} \left( z_{\nu}- y_{\nu}\right)^{\beta_{\nu} - 2\gamma_{\nu}}
\theta_{_{N}}(y)u(y) \, dy, 
\end{align*}
\begin{align*}
&e^{-\lambda \phi_{0}(w,z)}  I_{1,2}= \frac{1}{(2N)!}  e^{-\frac{\lambda}{2} \left[\left( w''\right)^{2} + \left( z''\right)^{2}\right]} 
\\
&
\times
\int e^{-\frac{\lambda}{2}(w-s)^{2}}  \chi_{0}(s)s^{2N} ds
\sum_{ |\beta| \leq 2N}\sum_{\substack{\gamma \leq \lfloor \frac{\beta}{2}\rfloor \\  \gamma \in \mathbb{Z}^{n}_{+} }}
\frac{\beta!\, i^{|\beta|-2|\gamma|}}{\gamma! \left( \beta-2\gamma!\right) 2^{\gamma}}
\lambda^{|\beta|-|\gamma|} 
\\
&
\times
\int_{ \mathscr{B}_{\tilde{r}_{0}}(x_{0})}
\!\!\!\left[ e^{-\frac{\lambda}{2}(z-\zeta)^{2}}
\tilde{a}_{_{2N,\beta}}\!\! \left( \zeta \right) 
\prod_{\nu=1}^{n}\! \left( z_{\nu}-\zeta_{\nu}\right)^{\beta_{\nu} - 2\gamma_{\nu}}
\!\tilde{u}(\zeta) \right]_{ \zeta = y + i\varepsilon_{4} \vartheta(y)\eta^{0}  }
\hspace{-4em}
\det\left(B(y,\varepsilon_{4},\eta^{0})\right)
\, dy,
\end{align*}
%
	where $ B(y,\varepsilon_{4},\eta^{0})$ is the $n\times n$-matrix 
%
\vspace{0.3em}
\begin{equation*}
\tiny{
	\begin{pmatrix}
	1 +i \varepsilon_{4}\, \vartheta^{(e_{1})} (y)\eta^{0}_{1} \!\! & \!\! i \varepsilon_{4} \vartheta^{(e_{2})} (y)\eta^{0}_{1}\!\!  & \!\!\! \cdots & \!\!\!  \cdots & \!\!\!  \cdots
	&\!\!\!  i \varepsilon_{4} \vartheta^{(e_{n})} (y)\eta^{0}_{1} 
	\\
	i \varepsilon_{4} \vartheta^{(e_{1})} (y)\eta^{0}_{1} \!\! & \!\!\! 1 +i \varepsilon_{4} \vartheta^{(e_{2})} (y)\eta^{0}_{2} &
	\!\! i \varepsilon_{4} \vartheta^{(e_{3})} (y)\eta^{0}_{2}\!\! &\!\!\!   \cdots \!\!\! &\!\!\!   \cdots \!\!\! &\!\!\!    i \varepsilon_{4}  \vartheta^{(e_{n})} (y)\eta^{0}_{2} \\
	\vdots  & \!\!\!  \ddots & \!\!\! \ddots & \!\!\! \ddots &\!\!\!  \ddots  &\!\!\!  \vdots  \\
	i \varepsilon_{4} \vartheta^{(e_{1})} (y)\eta^{0}_{n} \!\!\! & \!\!\!  \cdots\!\!\!  & \!\!\! \cdots\!\!\!  & \!\!\! \cdots\!\!\!  &\!\!\!  i \varepsilon_{4} \vartheta^{(e_{n-1})} (y)\eta^{0}_{n}
	&\!\!\!  1+  i \varepsilon_{4} \vartheta^{(e_{n})} (y)\eta^{0}_{n}
	
	\end{pmatrix}
}
\end{equation*}
\vspace{0.3em}
	and
%
\begin{align*}
e^{-\lambda \phi_{0}(w,z)}  I_{1,3}=
&
\frac{1}{(2N)!}  e^{-\frac{\lambda}{2} \left[\left( w''\right)^{2} + \left( z''\right)^{2}\right]} 
\int e^{-\frac{\lambda}{2}(w-s)^{2}}  \chi_{0}(s)s^{2N} ds
\\
&
\times
\sum_{ |\beta| \leq 2N}\sum_{\substack{\gamma \leq \lfloor \frac{\beta}{2}\rfloor \\  \gamma \in \mathbb{Z}^{n}_{+} }}
\frac{\beta!\, i^{|\beta|-2|\gamma|}}{\gamma! \left( \beta-2\gamma!\right) 2^{\gamma}}
\lambda^{|\beta|-|\gamma|} 
\sum_{\ell=1}^{n}  
I_{1,3,\beta,\gamma,\ell},
\end{align*}
%
	where 
%
\begin{align*}
&I_{1,3,\beta,\gamma,\ell}=
\\
& \int_{0}^{\varepsilon_{4}} \int_{ S_{\eta^{0},t} } 
\!\!\!\!\!\!\!\!\!\!  e^{-\frac{\lambda}{2}(z- \zeta)^{2}} \tilde{a}_{_{2N,\beta}} \!\! \left( \zeta \right) 
\prod_{\nu=1}^{n} \left( z_{\nu}- \zeta_{\nu} \right)^{\beta_{\nu} - 2\gamma_{\nu}}
\left( \bar{\partial}\tilde{u}\right) ( \zeta ) \,(i \vartheta\left( \frac{\zeta+\bar{\zeta}}{2}\right)\eta_{\ell}^{0}) \,dt \wedge d\zeta
\\
&=
\int_{0}^{\varepsilon_{4}} \int_{ \mathscr{B}_{\tilde{r}_{0}}(x_{0}) }
e^{-\frac{\lambda}{2}(z- y -it\vartheta(y)\eta^{0} )^{2}} \tilde{a}_{_{2N,\beta}} \!\! \left( y + it\vartheta(y)\eta^{0} \right) 
\left( \bar{\partial}\tilde{u}\right)( y + it\vartheta(y)\eta^{0} )
\\
&
\qquad\qquad\qquad
\times
\prod_{\nu=1}^{n} \left( z_{\nu}-   y_{\nu} -it\vartheta(y)\eta^{0}_{\nu} \right)^{\beta_{\nu} - 2\gamma_{\nu}}
\det\left(A_{\ell}(y,t,\eta^{0})\right)\,
dt\,dy.
\end{align*}
%
	Since $\left( z' - y \right)^{2} \geq \varepsilon_{0}$ in $\mathscr{B}_{r_{N}}(x_{0}) \setminus \mathscr{B}_{\tilde{r}_{0}}(x_{0})$,
	using the same strategy as for the term $I_{2}$, we obtain
%
\begin{align*}
\left| e^{-\lambda \phi_{0}(w,z)}  I_{1,1} \right| \leq C_{5} \left( \frac{1}{2}\right)^{N} e^{-\varepsilon_{2}\lambda} ,
\end{align*}
%
	where $C_{5}$ is a positive constant independent of $N$.\\
	A quick inspection of the terms $ I_{1,2}$ and $ I_{1,3}$ highlights that
	the main differences with respect to the already treated terms, are the behavior of the phase function on the integration path
	as well as the presence of $\bar{\partial}u$.
	We point out that setting $a^{\ell}_{p,q}(y,t,\eta^{0})$ and $b_{k,m}(y,\varepsilon_{4},\eta^{0})$, $p,q \in \left\{1,\dots, n+1 \right\}$
	and $k,m \in \left\{1,\dots, n \right\}$, the entries of the matrixes  $A_{\ell}(y,t,\eta^{0})$ and $B(y,\varepsilon_{4},\eta^{0})$ respectively,
	since we can estimate the entries $|a^{\ell}_{p,q}(y,t,\eta^{0})|$ and $|b_{k,m}(y,\varepsilon_{4},\eta^{0})|$ by $ (1 + \sup_{i}\|\vartheta^{(e_{i})}\|_{\infty})$
	we have 
%
\begin{align*}
&\left|\det \left(A_{\ell}(y,t,\eta^{0})\right) \right| \leq \sum_{\sigma \in S_{n+1}} \prod_{p=1}^{n+1}\left|a^{\ell}_{p,\sigma(q)}(y,t,\eta^{0})\right|
\\
&
\hspace{12em}
\leq \frac{(n+1)!\left[(n+1)!+1\right]}{2}  (1 + \sup_{i}\|\vartheta^{(e_{i})}\|_{\infty})^{n+1},
\\
&
\left|\det \left(B(y,\varepsilon_{4},\eta^{0})\right) \right| \leq \sum_{\sigma \in S_{n}} \prod_{k=1}^{n}\left|b_{k,\sigma(k)}(y,\varepsilon_{4},\eta^{0})\right|
\\
&
\hspace{17em}
\leq \frac{n!\left(n!+1\right)}{2}  (1 + \sup_{i}\|\vartheta^{(e_{i})}\|_{\infty})^{n}.
\end{align*}
%
	We focus on the exponential function:
%
\begin{align*}
e^{-\frac{\lambda}{2}\Re(z-  y - i t \vartheta(y)\eta^{0} )^{2}} 
=
e^{\frac{\lambda}{2} \left( z''\right)^{2} } 
e^{ - \frac{\lambda}{2} \left( z' - y \right)^{2} } 
e^{- \lambda t \vartheta(y) z''\eta^{0} + \frac{\lambda}{2} \left( t \vartheta(y)\right)^{2} |\eta^{0}|^{2} }.
\end{align*}
%
	Since $z''$ is in a neighborhood of $-\xi_{0}$ then $z''\eta^{0} >0$. Hence there is a positive constant $c$ such that $ z''\eta^{0} > c |z''| |\eta^{0}|$;
	moreover since we can assume that there is a strictly positive constant $ a $ such that $|z''| \geq a$ then  $ z''\eta^{0} > c_{1} |\eta^{0}|$, $c_{1} >0$.
	We can estimate the above quantity with
%
\begin{align*}
e^{\frac{\lambda}{2} \left( z''\right)^{2} } 
e^{ - \frac{\lambda}{2} \left[ \left( z' - y \right)^{2}  
	+  	t \vartheta(y) |\eta^{0}| \left( 2c_{1} -  t \vartheta(y) |\eta^{0}| \right) \right]}.
\end{align*} 
%
	Choosing $t$ sufficiently small we have that $ 2c_{1} -  t \vartheta(y) |\eta^{0}| >0$.
	In the case $t=\varepsilon_{4}$,
	we obtain the analytic exponential decay for $I_{1,2}$; more precisely the same strategy used to handle the term $I_{2}$
	gives that there are two positive constants $C_{6}$ and $\tilde{\varepsilon}_{2}$, independent of $N$, such that
%
\begin{align*}
\left| e^{-\lambda \phi_{0}(w,z)}  I_{1,2} \right| \leq C_{6} \left( \frac{1}{2}\right)^{N} e^{-\tilde{\varepsilon}_{2}\lambda}.
\end{align*}
%
	In order to estimate the last term, $ \left| e^{-\lambda \phi_{0}(w,z)}  I_{1,3} \right| $, we can apply once again the strategy used to estimate $I_{2}$.
	The only difference is that we have to take care of the term $\left| \left( \bar{\partial}\tilde{u}\right)(y+it\vartheta(y)\eta^{0}) \right| $.
	Keeping in mind that $\tilde{u}$ is an $(s_{0}-1)$-almost analytic extension of $u$, we have
%
\begin{align*}
&\left| e^{-\frac{\lambda}{2}(z-\zeta)^{2}}\right|  \left| \left( \bar{\partial}\tilde{u}\right)(y+it\vartheta(y)\eta^{0}) \right|
\\
&
\qquad
\leq C
e^{\frac{\lambda}{2} \left( z''\right)^{2} } 
e^{ - \frac{\lambda}{2} \left( z' - y \right)^{2} } 
e^{- \lambda c_{2} t \vartheta(y)|\eta^{0} |} 
e^{- \varepsilon_{_{K}} \left( t \vartheta(y) |\eta^{0}|\right)^{-\frac{1}{s_{0}-1}}}
\\
&
\qquad
\leq C e^{\frac{\lambda}{2} \left( z''\right)^{2} - \frac{\lambda}{2} \left( z' - y \right)^{2} } e^{- \tilde{\varepsilon}_{_{K}}\lambda^{1/s_{0}}},  
\end{align*}
%
	where $ \tilde{\varepsilon}_{_{K}} = c_{2} \gamma_{1}^{(s_{0}-1)/s_{0}} + \gamma_{1}^{-1/s_{0}}$, $\gamma_{1}  = \varepsilon_{K}/(c_{2}(s_{0}-1))$,
	and $\varepsilon_{_{K}}$ is as in the Definition \ref{Def-almost} with $K = \overline{\mathscr{B}_{\tilde{r}_{0}}(x_{0})}$.
	The estimate in the exponential is obtained taking $\displaystyle\inf_{b} \left( \lambda c_{1} b + \varepsilon_{_{K}} b^{-\frac{1}{s_{0}-1}}\right)$,
	where $b = t\vartheta(y)|\eta^{0}|$.
	Using this estimate we conclude that there are two positive constants $C_{7}$ and $\varepsilon_{4}$ such that
%
\begin{align*}
\left| e^{-\lambda \phi_{0}(w,z)}  I_{1,3} \right| \leq C_{7} \left( \frac{1}{2}\right)^{N} e^{- \varepsilon_{4}\lambda^{1/s_{0}}}.
\end{align*}
%
	We deduce that there is a positive constant $C_{8}$ such that
%
\begin{align}\label{Est_I_1}
|e^{-\lambda \phi_{0}(w,z)} I_{1}|
\leq C_{8} \left( \frac{1}{2}\right)^{N} e^{- \varepsilon_{4}\lambda^{1/s_{0}}}.
\end{align}
\begin{remark}\label{RK_def_2}
	The estimate of the second term on the right hand side of (\ref{Split_int}), i.e. in the region $\mathscr{B}_{\tilde{r}_{0}}(x_{0})$,
	can be obtained in a similar way introducing  the family of homeomorphisms
	$$ \mathscr{H}_{t} : \mathscr{B}_{\tilde{r}_{0}}(x_{0}) \ni y \rightarrow \left( y_{1} + i t\eta^{0}_{1}, \dots,  y_{n} + i t\eta^{0}_{n}\right) \in \mathbb{C}^{n}_{\zeta},$$
	where $\eta^{0} \in \Gamma_{\varepsilon_{2}}$. Also in this case 
	$ \mathscr{H}_{t} \left( \mathscr{B}_{\tilde{r}_{0}}(x_{0}) \right)$ is a $n$-dimensional manifold of $\mathbb{C}^{n}_{\zeta}$ for every $t \in [0,1]$.
	Setting
	\begin{align*}
	V(y, t\eta^{0}) = \tilde{a}_{_{2N,\beta}} (y+it\eta^{0})  &e^{-\frac{\lambda}{2}(z-y-it\eta^{0})^{2}}\times \\
	&
	\prod_{\nu=1}^{n} \left( z_{\nu}- y_{\nu}-it\eta^{0}_{\nu}\right)^{\beta_{\nu} - 2\gamma_{\nu}} \tilde{u}(y+it\eta^{0}),
	\end{align*}
	and applying, also in this case, the Stokes' theorem
	\begin{align*}
	\int_{\mathscr{H}_{1} \left( \mathscr{B}_{\tilde{r}_{0}}(x_{0}) \right)} V(\zeta) \, d\zeta
	-
	\int_{\mathscr{H}_{0} \left( \mathscr{B}_{\tilde{r}_{0}}(x_{0}) \right)} V(\zeta) \, d\zeta
	=
	\int_{\mathscr{V}} d\left(V(\zeta)\right) \wedge d\zeta,
	\end{align*}
	where $\mathscr{V} = \left[0,1\right] \times \mathscr{H}_{0} \left( \mathscr{B}_{\tilde{r}_{0}}(x_{0}) \right)$,
	the estimate (\ref{Est_I_1}) can be obtained following step by step the strategy employed above.
\end{remark}

\vspace*{1em}
\noindent
	By (\ref{Est-I_2}), (\ref{Est-I_3}), (\ref{Est-I_4}), (\ref{Est-I_5}) and (\ref{Est_I_1}) we have
%
\begin{align*}
&|e^{-\lambda \phi_{0}(w,z)} \mathscr{P}_{N}(w,z)|
\leq 
C_{8} \left( \frac{1}{2}\right)^{N} e^{- \varepsilon_{4}\lambda^{1/s_{0}}}
+ C_{2} \left( \frac{1}{2}\right)^{N} e^{-\varepsilon_{2}\lambda}
\\
&
\hspace{10em}
+ \left[C_{3} (N+1)
+ C_{4} 
+C_{5}  \right] (N+1) \left( \frac{1}{2}\right)^{N} e^{-\varepsilon_{2}\lambda}.
\end{align*}

\noindent
	Summing up we obtain that there are two positive constants $C$ and $\varepsilon$ such that
%
\begin{align*}
| e^{-\lambda \phi_{0}(w,z)} T\left(\chi U\right)(w,z, \lambda) | \leq C e^{- \varepsilon \lambda^{1/s_{0}}},
\end{align*}
%
	for all $(w,z)$ in a neighborhood of $(0, x_{0}-i\xi_{0}) \in \mathbb{C}^{1+n}$. 
%
%
\vspace{1.5em}

\noindent
	\textbf{Step two}: if  $(0,x_{0},0,\xi_{0}) \notin WF_{s_{0}}(U)$ then $(x_{0},\xi_{0})\notin WF_{s_{0}}(u)$.
	In the analytic category the result was obtained in \cite{bcm_1990} via Fourier transform and taking advantage from the Theorem 8.2.4 in \cite{H_Book-1}. 
	Via FBI transform it is a consequence of a result in \cite{Iag_75} on the restriction of a distribution to a sub-manifold.
	More precisely we remark that for every  $\tau_{0} \neq 0$ the points of the form $(t_{0},x_{0},\tau_{0},\xi_{0})$ do not belong to $WF_{s_{0}}(U)$
	for every $s_{0}\geq 1$. This can be obtained ether via FBI transform, performing the classical deformation argument of the
	integral path with respect to the $t$-variable, or noticing that the operator $Q$ is elliptic for $\tau\neq 0$.
	Since  $ WF_{s_{0}}\left( \delta(t)\right) = \lbrace (x,0,0,\tau): \, x\in \mathbb{R}^{n} \text{ and } \tau \in \mathbb{R}\setminus\lbrace 0 \rbrace  \rbrace$
	we have that $ WF_{s_{0}}(U) \cap  WF_{s_{0}}\left( \delta(t)\right) = \emptyset $, or equivalently that the normal to the manifold $ t = 0 $ does not intersect 
	the $WF_{s_{0}} (U)$,  then the product of $U$ and $ \delta(t)$ is well defined.
	This allow us to consider $u(x)$ as $U(t,x) \times \delta(t)$ in the sense of distributions. 
	More in general we can define the map $\pi : \lbrace U \in \mathscr{E}'( \mathbb{R}^{n+1}) : \, WF_{s_{0}}(U) \cap  WF_{s_{0}}\left( \delta(t)\right) = \emptyset \rbrace
	\rightarrow \mathscr{E}'( \mathbb{R}^{n}) $ in the following way $u(\phi_{0})=\pi(U)(\phi_{0}) = U( \phi_{1} \delta(t)) $
	for all $\phi_{0} \in C_{0}^{\infty}(\mathbb{R}^{n})$ where $\phi_{1} \in C_{0}^{\infty}(\mathbb{R}^{n+1})$ and $\pi ( \phi_{1}) = \phi_{0}$. 
	Following the same strategy used in \cite{H_Book-1} we have that
	$WF_{s_{0}}( u) = WF_{s_{0}}(\pi(U))$ which is contained in
	$\lbrace (x,\xi) \in \mathbb{R}^{n}\times \mathbb{R}^{n}\setminus \lbrace 0 \rbrace \,:\, \exists \tau \in \mathbb{R}  \text{ with } (x,0,\xi,\tau) \in WF_{s_{0}}(U)\rbrace$.\\

\noindent
	This concludes the proof of Theorem \ref{Miro_AnV}.
%
%
\bibliographystyle{cdraifplain}

\begin{thebibliography}{00}
	
	
	%
	\bibitem{ABC}
	{\sc P.~Albano, A.~Bove, G.~Chinni},
	{\it Minimal Microlocal Gevrey Regularity for ``Sums of Squares'' },
	Int. Math. Res. Notices, {\bf 12} (2009), 2275-2302.
	%
	\bibitem{BG-72}
	{\sc M. S. Baouendi and C. Goulaouic},
	\emph{Nonanalytic-hypoellipticity for some degenerate elliptic operators},
	Bull. Amer. Math. Soc., \textbf{78} (1972), 483--486.
	%
	\bibitem{BCH_book} 
	{\sc S. Berhanu, P. D. Cordaro and J. Hounie},
	\emph{ An Introduction to Involutive Structures},
	Cambridge University Press, 2008.
	%
	\bibitem{Ber_Hai_2017}
	{\sc S. Berhanu and A. Hailu},
	\emph{Characterization of Gevrey regularity by a class of FBI transforms},
	In: Pesenson I., Le Gia Q., Mayeli A., Mhaskar H., Zhou DX. (eds) Recent Applications of Harmonic Analysis to Function Spaces, Differential Equations, and Data Science; Applied and Numerical Harmonic Analysis; Birkh\"auser, Cham (2017), 451--482.
	%
	\bibitem{BH_2012}
	{\sc S. Berhanu and J. Hounie},
	\emph{A Class of FBI Transforms}
	Commun. Partial Diff. Eq. {\bf 37}(1) (2012), 38--57.
	%
	\bibitem{bc_1981}
	{\sc P. Bolley and J. Camus},
	{\it Regularite Gevrey et iteres pour une classe d'operatours hypoelliptiques}
	Commun. Partial Diff. Eq. {\bf 6}(10) (1981), 1057--1110.
	%
	\bibitem{bcm_1990}
	{\sc P. Bolley, J. Camus and G. Metivier},
	{\it Th\'eor\`eme d'unicit\'e pour des vecteurs analytiques}
	Journal of differential equations, {\bf 86}(1)(1990), 59--72. 
	%
	\bibitem{BCN-82} {\sc P. Bolley, J. Camus and J. Nourrigat},
	{\it La condition de H\"ormander-Kohn pour les op\'erateurs pseudo-diff\'erentiels},
	Commun. Partial Diff. Eq. {\bf 7} (1982), 197-221.
	%
	\bibitem{bccj_2016}
	{\sc N. Braun Rodrigues, G. Chinni, P. D. Cordaro, and M. R. Jahnke},
	{\it Lower order perturbation and global analytic vectors for a class of globally analytic hypoelliptic operators},
	Proc. Amer. Math. Soc. {\bf 144} (2016), no. 12, 5159--5170.
	%
	\bibitem{Ch_PHD_t}
	{\sc G. Chinni}
	{\it Analytic and Gevrey (micro-)Hypoellipticity for Sums of Squares: an FBI Approach},
	Phd Thesis, 2008.
	%
	\bibitem{DH-1980}
	{\sc M. Damlakhi and B. Helffer},
	{\it Analyticit\'e et it\`eres d'un syst\`eme de champs non elliptique},
	Ann. scient. \'Ec. Norm. Sup.,
	$4^{e}$ s{\'e}rie, {\bf 13} (1980), 397--403.
	%
	\bibitem{DZ_1973}
	{\sc M. Derridj and C. Zuily}, 
	{\it R\'egularit\'e analytique et Gevrey d'op\'erateurs elliptiques d\'eg\'en\'er\'es},
	J. Math. Pures Appl. {\bf 52} (1973), 309-336.
	%
	\bibitem{D_2018}
	{\sc M. Derridj},
	{\it Local estimates for H\"ormander's operators with Gevrey coefficients and application to the regularity of their Gevrey vectors},
	Tunisian Journal of Mathematics, Vol. 1, No. 3 (2019), 321--345.
	%
	\bibitem{D_2019}
	{\sc M. Derridj}
	{\it Local estimates for H\"ormander's operators of first kind with analytic Gevrey coefficients and application to the regularity of their Gevrey vectors},
	Pacific Journal of Mathematics  {\bf 302} (2019), No. 2, 511--543.
	%
	\bibitem{GS} 
	{\sc A. Grigis and J. Sj\"ostrand}, 
	{\it Front d'onde analytique et somme de carr\`es de champs de vecteurs\/}, 
	Duke Math. J. {\bf 52} (1985), 35-51.
	%
	\bibitem{HM_1980}
	{\sc B. Helffer and CI. Mattera},
	{\it Analyticite et iteres reduits d'un systeme de champs de vecteurs},
	Communications in Partial Differential Equations, Vol. 5, No. 10 (1980), 1065--1072.
	%
	\bibitem{H67}
	{\sc L. H\"ormander}, 
	{\it Hypoelliptic second order differential equations}, 
	Acta Math. {\bf 119} (1967), 147-171.
	%
	\bibitem{H_FIO}
	{\sc L. H\"ormander}, 
	{\it Fourier integral operators. I}, 
	Acta Math. {\bf 127} (1971), 79--183.
	%
	\bibitem{H_Book-1}
	{\sc L. H\"ormander}
	{\it The Analysis of Linear Partial Differential Operators I},
	Springer-Verlag Berlin Heidelberg (1983).
	%
	\bibitem{Iag_75}
	{\sc D. Iagolnitzer},
	{\it Appendix Microlocal essential support of a distribution and decomposition theorems -- An introduction},
	In: Pham F. (eds) Hyperfunctions and Theoretical Physics. Lecture Notes in Mathematics, vol 449. Springer, Berlin, Heidelberg (1975).
	%
	\bibitem{RS}
	{\sc L. Preiss Rothschild and E. M. Stein}, 
	{\it Hypoelliptic differential operators and nilpotent groups\/},
	Acta Math.  {\bf 137} (1976), 247-320.
	%
	\bibitem{Sj-Ast} 
	J. Sj\"ostrand, 
	{\it Singularit\'es analytiques microlocales\/}, Ast\'erisque {\bf 95} (1982).
	%
	
	
\end{thebibliography}
%

%
\end{document}